\pdfoutput=1
\RequirePackage{ifpdf}
\ifpdf 
\documentclass[pdftex]{sigma}
\else
\documentclass{sigma}
\fi

\numberwithin{equation}{section}

\newtheorem{Theorem}{Theorem}[section]
\newtheorem{Corollary}[Theorem]{Corollary}
\newtheorem{Lemma}[Theorem]{Lemma}
\newtheorem{Proposition}[Theorem]{Proposition}
 { \theoremstyle{definition}
\newtheorem{Definition}[Theorem]{Definition}

\newtheorem{Example}[Theorem]{Example}
\newtheorem{Remark}[Theorem]{Remark} }

\def\Sym{\mathrm{Sym}}
\def\sym{\mathfrak{s}\mathrm{y}\mathfrak{m}}
\usepackage[all]{xy}

\begin{document}
\allowdisplaybreaks

\renewcommand{\thefootnote}{}

\newcommand{\arXivNumber}{2101.07470}

\renewcommand{\PaperNumber}{016}

\FirstPageHeading

\ShortArticleName{Darboux Transformations for Orthogonal Differential Systems}

\ArticleName{Darboux Transformations for Orthogonal Differential\\ Systems and Differential Galois Theory\footnote{This paper is a~contribution to the Special Issue on Symmetry, Invariants, and their Applications in honor of Peter J.~Olver. The~full collection is available at \href{https://www.emis.de/journals/SIGMA/Olver.html}{https://www.emis.de/journals/SIGMA/Olver.html}}}

\Author{Primitivo ACOSTA-HUM\'ANEZ~$^{\rm a}$, Moulay BARKATOU~$^{\rm b}$, Raquel S\'ANCHEZ-CAUCE~$^{\rm c}$\newline and Jacques-Arthur WEIL~$^{\rm b}$}

\AuthorNameForHeading{P.~Acosta-Hum\'anez, M.~Barkatou, R.~S\'anchez-Cauce and J.-A.~Weil}

\Address{$^{\rm a)}$~Instituto de Matem\'atica \& Escuela de Matem\'atica, Universidad Aut\'onoma de Santo Domingo,\\
\hphantom{$^{\rm a)}$}~Dominican Republic}
\EmailD{\href{mailto:pacosta-humanez@uasd.edu.do}{pacosta-humanez@uasd.edu.do}}

\Address{$^{\rm b)}$~XLim - Universit\'e de Limoges \& CNRS, Limoges, France}
\EmailD{\href{mailto:moulay.barkatou@unilim.fr}{moulay.barkatou@unilim.fr}, \href{mailto:jacques-arthur.weil@unilim.fr}{jacques-arthur.weil@unilim.fr}}

\Address{$^{\rm c)}$~Department of Artificial Intelligence, Universidad Nacional de Educaci\'on a Distancia\\
\hphantom{$^{\rm c)}$}(UNED),~Madrid, Spain}
\EmailD{\href{mailto:rsanchez@dia.uned.es}{rsanchez@dia.uned.es}}

\ArticleDates{Received July 21, 2022, in final form March 20, 2023; Published online March 31, 2023}

\Abstract{Darboux developed an ingenious algebraic mechanism to construct infinite chains of ``integrable'' second-order differential equations as well as their solutions. After a surprisingly long time, Darboux's results were rediscovered and applied in many frameworks, for instance in quantum mechanics (where they provide useful tools for supersymmetric quantum mechanics), in soliton theory, Lax pairs and many other fields involving hierarchies of equations. In this paper, we propose a method which allows us to generalize the Darboux transformations algorithmically for tensor product constructions on linear differential equations or systems. We obtain explicit Darboux transformations for third-order orthogonal systems ($\mathfrak{so}(3, C_K)$ systems) as well as a framework to extend Darboux transformations to any symmetric power of $\mathrm{SL}(2,\mathbb{C})$-systems. We introduce SUSY toy models for these tensor products, giving as an illustration the analysis of some shape invariant potentials.
All results in this paper have been implemented and tested in the computer algebra system \textsc{Maple}.}


\Classification{12H05; 35Q40; 81Q60}

\renewcommand{\thefootnote}{\arabic{footnote}}
\setcounter{footnote}{0}

\section{Introduction}
Consider a Sturm--Liouville problem
\begin{gather*}
(E_{\lambda}) \colon \ \frac{{\rm d}^2 y}{{\rm d}x^2}(x)+ q(x) y(x) = \lambda y(x).
\end{gather*}
Here, $ \lambda$ is a spectral parameter
and $q(x)$ is called a potential. Assume that, for $\lambda=0$, a non-zero solution $y_0(x)$ of $(E_0)$ is known.
Then, letting $\theta_0:=\frac{y_0'}{y_0}$, Darboux observed in 1882~\cite{da1}, see also~\cite{da,da2}, that the map $y\mapsto \tilde{y}:= y'-\theta_0 y$
maps solutions of $(E_{\lambda})$ to solutions of
\begin{gather*}
(\tilde{E}_{\lambda})\colon\ \frac{{\rm d}^2\tilde{y}}{{\rm d}x^2}(x) + \tilde{q}(x) \tilde{y}(x) = \lambda \tilde{y}(x),\qquad
	\textrm{with}\quad \tilde{q}(x) = q(x) + 2 \theta_0'.
\end{gather*}
We see that the transformation $({E}_{\lambda})\rightarrow \big(\tilde{E}_{\lambda}\big)$ is covariant: the potential $q$ is changed into $\tilde{q}$ but the shape of the equation is the same and the spectral parameter is untouched. Thanks to this covariance, one can iterate the transformation while preserving the shape of the equation. Thus, starting from a solvable case, for example, one obtains an infinite family of solvable equations.\footnote{In Theorem~\ref{tda1} below, we give this transformation for any second-order linear differential equation.}
To summarize, this transformation of Darboux may be viewed as a mapping on a couple formed of a generic solution and coefficients, here the potential $q$: it maps $(y,q)$ to $(\tilde{y},\tilde{q})$.

This seemingly harmless little transformation has been studied, applied and generalised by countlessly many authors and is now used to construct solutions and hierarchies of linear and nonlinear partial differential equations including the nonstationary Schr\"odinger equation, Korteweg--de Vries
and Kadomtsev--Petviashvili equations, $1 + 1$ and $2+1$ Toda lattice equations,
sine-Gordon and nonlinear Schr\"odinger equations: see the reference book of Matveev and Salle~\cite{MaSa91a} and also the books of
C.~Rogers and W.K.~Schief~\cite{RoSc02h} and C.~Gu, H.~Hu, and Z.~Zhou~\cite{GHZ}.
The common methodology is: start from a known solution and construct a covariant transformation to obtain infinite families of equations or solutions with prescribed properties and which all have the same shape.

The name \emph{Darboux transformation} seems to be strangely recent:~\cite{rosu1999short} and~\cite{GHZ} attribute it to Matveev in 1979 in the papers~\cite{Ma79u,Ma79m,MaSa79a}.
In these papers, the Darboux lemma concerning scalar Schr\"odinger
is extended to hierarchies
of linear partial differential equations of any orders (including the case of matrix or operator valued coefficients)
and their difference analogs, thus providing a systematic tool for construct explicit solutions
for hierarchies of linear and nonlinear PDE's and their difference versions.
In~\cite{Ma79a}, Matveev constructs multi-parametric families of real, singular and rational solutions of Zakharov--Schabat equations (for more general comments on this, see~\cite{Ma00a}).

Despite having been only recently revived, the Darboux transformation approach has had tremendous success and, nowadays, about 450 published papers in mathematics or mathematical physics have ``Darboux transformation'' in their title (not to mention at least 400 others which use Darboux transformations as a tool). Examples of such papers are the seminal works of Witten (supersymmetric quantum mechanics~\cite{wi}) and Gendenshte\"{\i}n (shape invariant potentials~\cite{ge}).
A matrix formalism for Darboux transformations for Schr\"odinger
equation was developed in the beginning of the twenty first century by Pecheritsin, Pupasov and Samsonov, see~\cite{sape}. Darboux transformations (and their matrix formalisms) also appear
in works of Novikov, Taimanov and Tsarev, see~\cite{NoTa16b,NoTa18v} and references therein.

Our contribution to this field takes its roots in a differential-Galoisian approach to Darboux transformations and shape
invariant potentials proposed by P.~Acosta-Humanez and \mbox{J.-J.~Morales}-Ruiz in~\cite{AP2,ac2,AP1}, where it was proved that the Darboux transformation
preserves the Galoisian structure of the differential equation (the Darboux transformation is isogaloisian). A~similar Galoisian approach was then explored by S.~Jim\'enez, \mbox{J.-J.~Morales}-Ruiz, R.~S\'anchez-Cauce and M.-A.~Zurro in~\cite{JiMoSaZu17r,JiMoSaZu18v,JiMoSaZu19f} in the context of integrable systems and rational solitons of KdV equations.
There, the authors studied the behavior of the Galoisian structure of some families of linear systems with respect to Darboux transformations.

One key to extending Darboux transformation to higher-order equations or systems is to view it as a Gauge transformation. For Sch\"odinger equations, this has been observed in~\cite{AP1} and~\cite{AP2,ac2}. For systems, such as AKNS systems, this approach can be found in the reference book of Gu, Hu and Zhu~\cite[Section~1.3, p.~18]{GHZ} and is used in papers such as~\cite{YaRu13m} or~\cite{JiMoSaZu17r,JiMoSaZu18v,JiMoSaZu19f}.
The contribution of this paper is to show how, viewing the Darboux transformation as a gauge transformation, tools from the Galois theory of differential equations (the Tannakian approach) allow us to construct Darboux transformations for linear differential equations or systems of order higher than two when they arise as constructions on an $\mathfrak{sl}(2)$-system, a situation that occurs in physical systems.

In Proposition~\ref{prop1}, we give a matrix factorization of this gauge transformation which in turn gives a simpler expression for the solutions of the new equation given by a Darboux transformation.
 This allows us to extend Darboux transformations to systems whose Galois group is represented as a symmetric power of $\mathrm{SL}(2,\mathbb{C})$. We then extend Darboux transformation to
orthogonal\footnote{A system $Y'=AY$ is called orthogonal when its matrix $A$ is in $\mathfrak{so}(3, \mathbb{C})$, i.e., $A+A^{\rm T}=0$ and $\textrm{trace}(A)=0$.} systems, i.e., $\mathfrak{so}(3, \mathbb{C})$-systems, which allows us to construct infinite chains of integrable (in Galois sense) linear differential systems.
 An important feature is that this gives us formulas which are simple and easy to use (generalizing Darboux transformations carelessly can easily produce unreadable formulas). The formulas are summarized in the diagrams~\eqref{eq:diag-1} and~\eqref{eq:diag-2}. Finally, we show in Section~\ref{extension} how these formulas may be generalized to any third-order linear differential system with an orthogonal Galois group, i.e., a Galois group in ${\rm SO}(3, \mathbb{C})$
 (or, equivalently, which has a quadratic first integral)). All formulas and results in this paper have been implemented and tested in the computer algebra system \textsc{Maple}.

 As an application of this approach, we study in Section~\ref{applications} shape invariant potentials in supersymmetric quantum mechanics for Schr\"odinger equations and we apply this formalism to recover some results from other authors (Fedorov, Maciejewski, Przybylska and others) related to $\mathfrak{so}(3, \mathbb{C})$ systems, in particular, the rigid solid
 problem and Frenet--Serret formulas in Section~\ref{sec:applic-so3}.

\section{Darboux transformation and differential Galois theory}
In this section, we recast known material. We present the Darboux transformation for second-order equation; we then recall some basics of differential Galois theory and a Galoisian view of the Darboux transformation. We recall another classical result of Darboux on orthogonal systems and recall constructive methods to study higher-order representations of ${\rm SL}(2)$, namely symmetric powers. All this combined will allow us to construct Darboux transformations for families of orthogonal systems.

\subsection{Darboux transformation}
In~\cite{da1}, Darboux proposes a transformation which, given a family of linear differential equations produces a new family of differential equations with a similar shape and similar properties. This transformation has proved to be powerful, for example, in the study of Shr\"odinger equations.
We recast it here in modern language.
This proposition appears in Darboux's note~\cite{da1}.
Ince mentions it in~\cite[p.~132]{in}.

\begin{Theorem}[Darboux,~\cite{da1}]\label{tda1}
Consider the family of differential equations $\mathcal{L}(y)=m r y$:
\begin{gather}\label{etda11}
y''+ py'+(q -mr)y = 0,
\end{gather}
where $p$, $q$, $r$ are functions $(r\neq 0$ by hypothesis$)$ and $m$ is a constant parameter.
Given a~non-zero value for $m$, we let $y_m$ denote a general solution of~\eqref{etda11}.
Suppose that we know a~non-zero solution $y_0$ of equation~\eqref{etda11} for $m=0$.
Let $\tilde{y}_m$ be a function defined by 	
\begin{gather}\label{etda12}
\tilde{y}_m
 = \frac{1}{\sqrt{r}} (y_m'- \theta_0 y_m) \qquad \textrm{with} \quad
 \theta_0:= \frac{y_0'}{y_0}.
 \end{gather}
Then, for $m\neq 0$, $\tilde{y}_m$ is a general solution of the new differential equation
\begin{gather}\label{etda13}
\tilde{y}''+p\tilde{y}'+(\tilde{q}-mr)\tilde{y}=0 \qquad \textrm{with} \quad
\tilde{q} = q + q_0,
\end{gather}
where, letting $\hat{r}:=\frac{r'}{2r}$, the new part $q_0$ is given by
\begin{gather*}
q_0:= 2 \theta_{{0}}' + \hat{r}' + p'- \hat{r}(\hat{r} + 2 p -4\theta_0).
\end{gather*}
\end{Theorem}

The new $\tilde{q}$ is given by Darboux~\cite{da1} in the compact expression
\begin{gather*}
\tilde{q}=y_0 \sqrt{r}\bigg(\frac{p}{y_0 \sqrt{r}} -\bigg(\frac{1}{y_0 \sqrt{r}}\bigg)'\bigg)'.
\end{gather*}
This transformation has been made famous by its applications to Schr\"odinger equations where $p=0$ and $r=1$. In this case, the formula is much simpler: $\tilde{q}=q+ 2\theta_0'$.

\begin{Definition}
The map~\eqref{etda12} transforming the family of equations~\eqref{etda11} to the family~\eqref{etda13} is called the \emph{Darboux transformation}. 
\end{Definition}

As mentioned in the introduction, the philosophy of the Darboux transformation is to start from a known-solution to a Sturm--Liouville problem to construct a covariant transformation which yields a similar equation: the eigenvalue $m$ and the weight $r$ are preserved; only the potential $q$ is changed into $\tilde{q}$ and the shape of the equation. Because of this covariance, the transformation may be iterated, producing infinite families with common properties.

This transformation has been generalised by several authors and is now applied to construct solutions of linear and nonlinear partial differential equations including the nonstationary Schr\"odinger equation, Korteweg--de Vries and Kadomtsev--Petviashivili equations, $1 + 1$ and $2+1$ Toda lattice equations,
sine-Gordon and nonlinear Schr\"odinger equations: see the book of Matveev and Salle~\cite{MaSa91a}. The name \emph{Darboux transformation} seems to be rather recent:~\cite{rosu1999short} attributes it to Matveev in 1979 in the papers~\cite{Ma79u,Ma79m,MaSa79a}. Nevertheless, it has had tremendous success and, nowadays, more than 800 published papers in mathematics or mathematical physics have ``Darboux transformation'' in their title.

The approach developed here, though, seems to be new. We focus on families of differential equations which may be viewed as constructions from ${\rm SL}(2)$ and give a methodology to generate Darboux transformation for them. In order to build this, we first start with some differential Galois theory.

\subsection{Differential Galois groups} \label{dgal-section}
Differential Galois theory, also known as Picard--Vessiot theory, is analogous to the classical Galois theory for polynomials; it describes algebraic relations that may exist between solutions of linear differential equations and their derivatives,
see~\cite{crhamo,PS}.
A differential field $K$, depending on a variable $x$, is a field equipped with a derivation $\partial_x={}'$. We denote by $C_K$ the field of constants of $K$, the set of $c\in K$ such that $c'=0$. Along this
paper, we consider differential equations or systems whose coefficients belong
to a differential field $K$ whose constant field $C_K$ is algebraically closed and of characteristic zero.
For simplicity, we will explain this Galois theory on operators of order two but it applies similarly to linear differential equations or systems of any order.

Consider the differential operator
\begin{gather*}
\mathcal L:=\partial_x^2 +p \partial_x +q,\qquad p,q\in K.
\end{gather*}
Let $\{ y_1,y_2\}$ denote a basis of solutions of $\mathcal Ly=0$.
We let $F:=K(y_1,y_2,y_1',y_2')$ be the smallest differential extension of $K$
containing these solutions of $\mathcal Ly=0$ and such that $C_K=C_L$.
The differential extension $F$ is called a \emph{Picard--Vessiot extension} of $K$ associated to $\mathcal Ly=0$.
The differential automorphisms of $F$ are the automorphisms that commutes with the derivation. The differential Galois group of $\mathcal Ly=0$, denoted by $\mathrm{DGal}(F/K)$, is the group of differential automorphisms of $L$ which leave invariant each element of $K$.

Let $\sigma\in \mathrm{DGal}(F/K)$. Then, $\{\sigma(y_1),\sigma(y_2)\}$ is another basis of solutions of $\mathcal Ly=0$. It~follows that
$\sigma (y_1, y_2)=(y_1,y_2)M_\sigma$. This $M_\sigma$ is the matrix of the automorphism $\sigma$.
We see that $\mathrm{DGal}(F/K)$ is a group of matrices and it is actually a linear algebraic group,
\begin{gather*}
\mathrm{DGal}(F/K) 	\subseteq \mathrm{GL}(2, C_K).
\end{gather*}
The Wronskian $W$ of the solutions $y_1$ and $y_2$
satisfies the differential equation $W'+pW=0$. Thus, $W=\exp(\int (-p){\rm d}x)$. We find that $W\in K$ if and only if $p=\frac{w'}{w}=(\ln w)'$ for some $w\in K$. In this case, we have $\sigma (W)=W$ (because $W\in K$).
As $\sigma (W)=\det(M_\sigma)W$ with $W\neq 0$, we obtain $\det(M_\sigma)=1$, that is,
\begin{gather*}
\mathrm{DGal}(F/K)\subseteq \mathrm{SL}(2, C_K) \
\Longleftrightarrow \ p=\frac{w'}{w } = (\ln w)', \qquad w\in K.
\end{gather*}

We say that an algebraic group $G$ is \emph{virtually solvable} when the connected identity component of $G$, denoted by $G^\circ$, is a solvable group. In this paper, we say that $\mathcal Ly=0$ is \emph{integrable} whenever $\mathrm{DGal}(F/K)$ is virtually solvable. This corresponds to cases when one can compute formulas for the solutions.

Given a non-zero solution $y$ and some non-zero function $c$, the standard change of variables (see, e.g., the book of Ince~\cite{in})
\begin{gather*}
u=-\frac1{c} \cdot \frac{y'}{y}
\end{gather*}
changes our second-order linear differential equation to the first-order (nonlinear) Riccati equation
\begin{gather*}
(R)\colon \ u'=a+bu+cu^2, \qquad \textrm{where}\quad b = -p - \frac{c'}{c}\quad \textrm{and}\quad a=\frac1{c} q.
\end{gather*}
Differential Galois theory shows that the equation $\mathcal{L}(y)=0$ is integrable if and only if the Riccati equation $(R)$ has an algebraic solution. Similar statements (although more technical) can be obtained for higher-order equations, see~\cite{PS}.

A scalar differential equation such as $\mathcal Ly=0$ is equivalent to its companion linear differential system $[A]$:
\begin{gather*}
[A]\colon\ X'=-AX,\qquad \text{where}\quad A=\begin{pmatrix}0&-1\\ q& p\end{pmatrix}
\quad \text{and}\quad X=\begin{pmatrix} y \\ y'\end{pmatrix}\!.
\end{gather*}
Differential Galois theory applies naturally to linear differential systems as well (solution spaces are vector spaces, the groups act on these vector spaces).

Given a linear differential system $[A]\colon X'=-AX$, a \emph{gauge transformation} is a linear change of variables $X=PY$ with $P \in {\rm GL}(n,K)$. It transforms the system $[A]$ into a linear differential system $Y'=-P[A] Y$ with
\begin{gather*}
P[A]:= P^{-1} A P + P^{-1} P'.
\end{gather*}
We say that two linear differential systems $[A]$ and $[B]$ are \emph{equivalent} over $K$ when there exists a~gauge transformation $P \in {\rm GL}(n,K)$ such that $B=P[A]$. Equivalent differential systems share the same differential Galois group.

Given any linear differential system $[A]\colon X'=-AX$, the \emph{cyclic vector method} (see~\cite{Ba93a} for a~simple constructive process) allows to construct an equivalent differential system in companion form. Namely, let $f$ denote any linear combination of the $n$ components of $X$; differentiating~$n$ times using the relations $X'=-AX$, we find that $f, f', \ldots, f^{(n)}$ are $n+1$ elements in a vector space of dimension $n$ (generated by the $n$ components of $X$) and hence they are linearly dependent: this gives a differential operator for $f$, see~\cite{Ba93a} and references therein for more on this \emph{cyclic vector} process. So we may go from operator to system and vice-versa without altering the theory.

The Darboux transformation is initially a transformation on operators. We will recast it as a transformation on systems and this will allow us to use the machinery on systems to build higher-order Darboux transformations.

\subsection{Darboux transformations and Galois groups}

Consider the differential field $K$ and the family of differential operators
\begin{gather}\label{eqdarb}
 {\mathcal{L}}_m:=\partial_x^2 + p \partial_x +(q -mr), \qquad p,q,r\in K\quad \textrm{with}\quad
 r\neq 0.
\end{gather}
Here, $m$ is a constant parameter.
When $m$ is given any value, we let $F_m$ be a Picard--Vessiot extension of $K$ corresponding to the equation ${\mathcal{L}}_m y=0$.
Without loss of generality,
 we assume from now on that there exists $w\in K$ such that $p=w'/w = (\ln w)'$.
 Therefore, the differential Galois groups $\mathrm{DGal}(F_m/K)$ of ${\mathcal{L}}_my=0$ are subgroups of $\mathrm{SL}(2,C_K)$.

Following Acosta-Hum\'anez and co-authors in~\cite{AP1}, see also~\cite{AP2, ac2}, we denote by $\Lambda$ the set of values of $m$ for which ${\mathcal{L}}_m y=0$ is integrable (over $K$), the so-called \emph{algebraic spectrum}.
 If we have a family of parameters $m$
 for which the equations ${\mathcal{L}}_m$ are integrable over $K$
(including $m=0$, so as to satisfy hypotheses of Darboux's theorem, Theorem~\ref{tda1} above), the Darboux transformation will construct a new family with the same shape, with no new transcendent in its coefficient field and equivalent with the first one, hence
preserving the integrability properties as explained below.

After performing a Darboux transformation, our new family $\widetilde{\mathcal{L}}_m$ of differential equations has a new
coefficients field

$\widetilde{K}:= K(\theta_0)$ (with notations of Theorem~\ref{tda1}).
Note that the Darboux transformation itself is defined over the bigger field $\widetilde{K}\big(\sqrt{r}\big)$.

We say that
the Darboux transformation is
\emph{isogaloisian} when the differential Galois group is preserved (see~\cite{AP2,ac2,AP1}), i.e.,
	$\mathrm{DGal}\big(\widetilde{F}_m/{\widetilde K}\big) =\mathrm{DGal}(F_m/K)$;
 it is \emph{strong isogaloisian} when $\widetilde K=K$
 (i.e., when $\theta_0\in K$).
 Results on the isogaloisian character of the Darboux transformation appear in~\cite{AP1}, see also~\cite{AP2, ac2}; they will reappear in the next section as a consequence of the view of the Darboux transformation as a gauge transformation. When $\theta_0$ is algebraic over $K$, then the Darboux transformation is virtually isogaloisian; in this case, for any value $m$ such that ${\mathcal{L}}_m$ is integrable, its Darboux transformation $\widetilde{\mathcal{L}}_m$ is also integrable (see~\cite{AP1}).
We see that the Darboux transformation transforms a family of integrable equations into another integrable family with the same shape.

\subsection{Third-order orthogonal systems}
We now turn to a rather different result, also due to Darboux, which allows one to solve third-order orthogonal systems using only solutions of a first-order Riccati equation, see
\cite[Part~I, Chapter~II]{da3} for details and proofs.
A third-order orthogonal system is one of the form
\begin{gather*}
 \begin{pmatrix} \alpha \\ \beta \\ \gamma \end{pmatrix}' = \begin{pmatrix} 0 & h & -g \\ -h & 0 & f \\ g & -f & 0 \end{pmatrix} \cdot \begin{pmatrix} \alpha \\ \beta \\ \gamma \end{pmatrix}\!.
\end{gather*}
A simple calculation shows that $\alpha^2+\beta^2+\gamma^2$ is always constant for such a system.

\begin{Theorem}[{Darboux,~\cite[Chapter~II, pp.~30--31]{da3}}]\label{tso3dar}
Consider a differential system
\begin{gather}\label{etso3dar}
 \begin{pmatrix} \alpha \\ \beta \\ \gamma \end{pmatrix}' = \begin{pmatrix} 0 & h & -g \\ -h & 0 & f \\ g & -f & 0 \end{pmatrix} \cdot \begin{pmatrix} \alpha \\ \beta \\ \gamma \end{pmatrix}\!.
\end{gather}
A solution $(\alpha, \beta, \gamma)$ such that $\alpha^2+\beta^2+\gamma^2=1$
can be parameterized by
\begin{gather}\label{etso3dar2}
\alpha=\dfrac{1-uv}{u-v}, \qquad
\beta={\rm i} \dfrac{1+uv}{u-v} \qquad \text{and} \qquad \gamma= \dfrac{u+v}{u-v},
\end{gather}
where $u$ and $v$ are distinct solutions of the same Riccati equation
\begin{gather*}
\theta'= \omega_0+\mu\theta+{\omega_1}\theta^2, \qquad \textrm{with} \quad
 \omega_0= \dfrac{g-{\rm i} f}{2}, \quad \omega_1= \dfrac{g+{\rm i} f}{2},
 \quad \mu=-{\rm i} h.
\end{gather*}
Furthermore,
\begin{gather*}
u=\dfrac{\alpha+{\rm i}\beta}{ 1-\gamma} = \dfrac{1+\gamma
}{\alpha-{\rm i}\beta} \qquad \text{and} \qquad v=-\dfrac{1-\gamma}{\alpha-{\rm i}\beta} = -\dfrac{\alpha+{\rm i}\beta}{1+\gamma}.
\end{gather*}
\end{Theorem}

\begin{Remark}\label{rmk-darboux-thm2}
As we recalled, solutions of a first-order Riccati equation are logarithmic derivatives of solutions of a second-order linear differential equation. So, this result shows that one can solve these third-order orthogonal systems using solutions of second-order linear differential equations.
Namely, performing the change of variable $u=-\frac{1}{\omega_1}\frac{y'}{y}$, we see that
$y$ is a solution of
\begin{gather}\label{linear-eq-darboux}
y''+\bigg( \mu \omega_1 + \frac{\omega_1'}{\omega_1}\bigg) y' + \omega_1 \omega_0 y = 0.
\end{gather}
So, the parametrization given by~\eqref{etso3dar2} can be restated as
\begin{gather*}
\begin{pmatrix}
\alpha \\ \beta \\ \gamma \end{pmatrix}	= \frac{1}{w}
\begin{pmatrix} - \omega_1 &0&{\frac {1}{\omega_{{1}}}}
\vspace{1mm}\\ -{\rm i}\omega_{{1}} &0&{\frac {-{\rm i}}{ \omega_{{1}} }}\\ 0&
1&0\end{pmatrix}
\begin{pmatrix} y_1y_2 \\ y_1'y_2+y_1y_2' \\ y_1'y_2'
\end{pmatrix}\!,
\end{gather*}
where $y_1$, $y_2$ are linearly independent solutions of the second-order linear differential equation~\eqref{linear-eq-darboux} and $w:=y_1'y_2-y_1y_2'$ is their Wronskian.
This gives an explicit solution for an orthogonal system in terms of solutions of second-order equations. This will be further explored in the next section.
\end{Remark}

In the study of the rigid solid, see for example Fedorov et~al.~\cite{femapr, femapr0}, orthogonal systems such as~\eqref{etso3dar} are traditionally written in a more compact way, using the cross-product $\times$:
\begin{gather}\label{fedo}
 Z'=Z \times \Omega,\qquad Z=(\alpha,\beta,\gamma)^{\rm T}, \qquad \Omega=(f,g,h)^{\rm T}, \qquad f,g,h\in K.
\end{gather}

From now on, we work with either equation~\eqref{fedo} or equation~\eqref{etso3dar}: they are the same equation, although presented differently.
In what follows, the terminology \emph{orthogonal systems} will refer to systems of the form~\eqref{fedo} or~\eqref{etso3dar}.

\subsection{Tensor constructions, invariants and symmetric powers}\label{tensor-constructions}

In differential Galois theory, one classically translates linear algebra constructions on the solution space (tensor product, symmetric powers, etc.) into constructions on differential systems. This allows to measure properties on solutions by looking for rational function solutions in these tensor constructions (see~\cite[Chapters~3 and~4]{PS}). We review here the construction of symmetric powers for later use.

Let $V$ denote a vector space over $C_K$ of dimension $n$. We fix a basis $\cal{B}$ of $V$ and consider $g\in \mathrm{End}(V)$. Let
$M = (m_{i,j}) \in {\mathcal M}_{n}(C_K)$ denote the $n\times n$ matrix of the endomorphism $g$ in that basis $\cal{B}$.
We define a linear action of $g$
on the variables $X_{j}$ by $g(X_{j}):=\sum_{i=1}^n m_{i,j} X_{i}$ so $g$ acts on the indeterminates $X_i$ as if they were the basis $\cal{B}$.

Consider a homogeneous polynomial
$P\in K[X_{1},\ldots,X_{n}]_{m}$ of degree $m$.
We may identify the $m$-th symmetric power of $V$ with the linear span of all monomials
$\{X_{1}^{m}, X_{1}^{m-1}X_{2}, \ldots, X_{n}^{m}\}$ of degree $m$. This way, our polynomial $P$ may be identified with its vector $v_{P}$ of coefficients on the monomial basis $\{X_{1}^{m}, X_{1}^{m-1}X_{2}, \ldots,X_{n}^{m}\}$.
Using the relations $g(X_{j}):=\sum_{i=1}^n m_{i,j} X_{i}$, we can define an action of
$g$ on $P$ by linear substitution $g(P):=P(g(X_{1}),\ldots, g(X_{n}))$.
This action by substitution translates into a natural action of $g$ on the coefficient vector $v_{P}$
by $g(v_{P}):=v_{g(P)}$.
Using this action, we define the \emph{$m$-th symmetric power in the sense of groups} $\Sym^m(M)$ of $M$ as the matrix of the linear map $v_{P} \mapsto v_{g(P)}$: it is defined by the relation
\begin{gather*}
v_{g(P)}=\Sym^m(M) \cdot v_{P}.
\end{gather*}
The map $M \mapsto \Sym^m(M)$ is a group morphism.
Given a group $G\subset {\rm GL}(V)$, an invariant of $G$ (in $\Sym(V)$) is a polynomial $P$ such that $\forall g\in G$, $g(P)=P$.

Similarly, one can associate to the matrix $M$ the derivation
\begin{gather*}
D_{M}=\sum_{j=1}^n \bigg(\sum_{i=1}^n m_{i,j} X_{i}\bigg)\frac{\partial}{\partial X_{j}}.
\end{gather*}
Then, we define the \emph{$m$-th symmetric power in the sense of Lie algebras} $\sym^{m}(M)$ of $M$ as the matrix of the linear map $v_{P} \mapsto v_{D_{M}(P)}$: it is defined by the relation
\begin{gather*}
v_{D_{M}(P)}=\sym^{m}(M) \cdot v_{P}.
\end{gather*}
The map $M \mapsto \sym^{m}(M)$ is a Lie algebra morphism.
Given a Lie algebra $\mathfrak{g}$, an invariant of $\mathfrak{g}$ is a polynomial
$P$ such that, $\forall D\in \mathfrak{g}$, $D(P)=0$. We have the equivalence: $P$ is an invariant of~$G^{\circ}$ if and only if $P$ is an invariant of its Lie algebra ${\rm Lie}(G)$.

Take a linear differential system $[A]\colon X' = -A X$. Its \emph{$m$-th symmetric power system}
is $[\sym^m (A)]$. If $X$ is a solution matrix of $[A]$, then $\Sym^m(X)$ is a solution matrix of
$[\sym^m (A)]$: $\Sym^m(X)'=-\sym^m(A)\Sym^m(X)$, see, e.g.,
\cite{ApCoWe13a, PS}.

\begin{Example}\label{ex-second-sym}
We will make this explicit on a construction used later on in the paper.
Consider a system $X'=-A X$ with
\begin{gather*}
A =\begin{pmatrix}0&-1 \\q&p\end{pmatrix}\!,\qquad \textrm{where} \quad p=\frac{w'}{w}
	\quad\textrm{with}\quad w\in K.
\end{gather*}
It is in companion form so it admits a fundamental solution matrix of the form
	\begin{gather*}
	 \mathbf{X}= \begin{pmatrix} y_1 & y_2 \\ y_1' & y_2' \end{pmatrix}\!.
	\end{gather*}
The second symmetric power system is $\big[ \sym^2 (A)\big]\colon Y'=-\sym^2 (A) Y$ with
\begin{gather*}
\sym^2 (A) = \begin{pmatrix} 0&-1&0\\2q &p
	&-2\\0&q& 2p
	\end{pmatrix}\!.
\end{gather*}
It admits the fundamental solution matrix
\begin{gather}\label{eq:mat ftal Y}
 \mathbf{Y}= \Sym^2 ( \mathbf{X}) = \begin{pmatrix} y_1^2 & y_1 y_2 & y_2^2 \\[3pt] 2y_1 y_1' & y_1' y_2 + y_1 y_2' & 2y_2y_2' \\[3pt] (y_1')^2 & y_1' y_2' & (y_2')^2 \end{pmatrix}\!.
\end{gather}
Now, let $\sigma \in \mathrm{DGal}(F_0/K)$ be an automorphism of the differential Galois group of $[A]$ with matrix representation
\begin{gather*}
M_{\sigma} = \begin{pmatrix} \lambda_{11} & \lambda_{12} \\ \lambda_{21} & \lambda_{22} \end{pmatrix} \in \mathrm{DGal}(F_0/K) , \qquad \lambda_{11} \lambda_{22} - \lambda_{12} \lambda_{21} = 1.
\end{gather*}
This means that $\sigma( \mathbf{X}) = \mathbf{X}\cdot M_{\sigma}$.
As $\Sym(\bullet)$ is a group morphism, we have
\begin{gather*}
\sigma( \mathbf{Y} ) = \sigma\big(\Sym^2( \mathbf{X})\big)
	= \Sym^2( \sigma( \mathbf{X}))
	= \Sym^2( \mathbf{X}\cdot M_{\sigma} )
	= \mathbf{Y}\cdot \Sym^2( M_{\sigma} ).
\end{gather*}
So the matrix of $\sigma$ acting on the solution $ \mathbf{Y}$ of $\big[\sym^2(A)\big]$ is
$\Sym^2( M_{\sigma})$ (computed as in formula~\eqref{eq:mat ftal Y}).
By a slight abuse of notation, we will denote the set of such matrices by $\Sym^2 (\mathrm{DGal}(F_0/K)) \subseteq {\rm SL}(3, C_K)$.
\end{Example}
Note that, using the cyclic vector $(1,0,0)^{\rm T}$ (see end of Section~\ref{dgal-section}), this symmetric power system can be written as the traditional third-order differential operator known as the second symmetric power of $\mathcal{L}$:
\begin{gather*}
\textrm{sym}^2 ( \mathcal{L} ):= \partial_x^3 +3p \partial_x^2 + \big(4q + p' +2p^2\big) \partial_x +2 (q' +2pq).
\end{gather*}

Last, we look at gauge transformations $P$.
The change $X = P Y$ transforms the system $[A]$ into $P[A]$. Then,
as $\Sym^m(\bullet)$ is a group morphism, $\Sym^{m}(X)=\Sym^{m}(P) \Sym^m(Y)$
which shows that $\sym^m (P[A]) = \Sym^{m}(P)[ \sym^m (A)]$.

\section{Darboux transformations for third-order orthogonal systems}

\subsection{Matrix formalism, Darboux transformation as a Gauge transformation}
A key to extending Darboux transformation to higher-order equations or systems is to view it as a Gauge transformation. For Sch\"odinger equations, this has been observed in~\cite{AP2,ac2,AP1} in a~Galoisian approach. For systems, such as AKNS systems, this approach can be found in the reference book of Gu, Hu and Zhu~\cite[Section~1.3, p.~18]{GHZ} and in papers such as~\cite{JiMoSaZu17r,JiMoSaZu18v,JiMoSaZu19f,YaRu13m}.
We now review this observation in a way that will allow us to construct our generalizations of Darboux transformations.

We write the differential equation ${\mathcal{L}}_my=0$, equation~\eqref{eqdarb}, in its companion system form:
\begin{gather}\label{mdar1}
[A_m]\colon\ X'=-(A_0+mN)X,\qquad\! X=\begin{pmatrix}y \\ y'\end{pmatrix}\!,\qquad\! A_0=\begin{pmatrix}0&-1 \\q&p\end{pmatrix}\!,\qquad\! N=\begin{pmatrix}0&0 \\ -r&0\end{pmatrix}\!.\!\!\!
\end{gather}
Note that $N^2$ is the null matrix.
The Darboux transformation transforms the family ${\mathcal{L}}_my=0$ into a family
$\widetilde{\mathcal L}_my=0$ whose companion form is now
\begin{gather}\label{mdar2}
\big[\widetilde A_m\big]\colon\ \widetilde{X}'=-\big(\widetilde A_0+mN\big) \widetilde{X}, \qquad
\widetilde{X} =\begin{pmatrix}\tilde{y} \\ \tilde{y}' \end{pmatrix}\!, \qquad \widetilde{A}_0=\begin{pmatrix}0&-1 \\\tilde{q}&p\end{pmatrix}\!,
\end{gather}
where $\tilde{y}$ and $\tilde{q}$ have the explicit form expressed in Theorem~\ref{tda1}. We denote by $\widetilde{K} = K(\theta_0)$ the field of coefficients of system~\eqref{mdar2}.

\begin{Proposition} \label{prop1}
The Darboux transformation given in Theorem~$\ref{tda1}$ is
equivalent to a
gauge transformation between the above families of systems $[\tilde{A}_m]$ and $[A_m]$,
whose matrix $P_m$ is given~by
\begin{gather}\label{gauge}
P_m:=\frac{1}{\sqrt{r}} \begin{pmatrix}
- \theta_0&1 \\ mr - \theta_0 \rho &\rho
\end{pmatrix}\!.
\end{gather}
The matrix $P_m$ factors as
\begin{gather*}
P_m= L_m . R, \qquad\textrm{where} \quad
L_m:= \begin{pmatrix}
0 & 1 \\ m r & \rho
\end{pmatrix}
\quad\textrm{and}\quad
R:=\frac{1}{\sqrt{r}}\begin{pmatrix}
1 & 0 \\ -\theta_0 & 1
\end{pmatrix}\!,
\end{gather*}
with $\theta_0=\frac{y_0'}{y_0}$, $p=\frac{w'}{w}$ and $\rho:=-\theta_0-p-\frac{r'}{2r}$.
When $\theta_0$ is algebraic over $K$, this gauge transformation preserves the identity component of the differential Galois group. 	
Moreover, whenever
$\theta_0 \in K$, this Darboux transformation
preserves the differential Galois group. 	
\end{Proposition}

\begin{Remark}
The new result, in this proposition, is the matrix factorization: the view as a gauge transformation and the isogaloisian properties come from~\cite{AP1}.
This factorization will be extended to our generalizations in the paper; it will allow us to obtain relatively simple formulas for our generalized Darboux transformations (general formulas can be cumbersome and unreadable otherwise).
Note that, in this factorization, the first matrix $L_m$ contains the dependence on $m$ and the second matrix $R$ only depends on the known solution $y_0$.
\end{Remark}

\begin{proof}
Let $Y:=P_m X$ with $Y=(z_1,z_2)^{\rm T}$ and $X$ given in the formula~\eqref{mdar1}. The first line is
$z_1=\frac{1}{\sqrt{r}}(-\theta_0 y + y')$ so we recognize the Darboux transformation from Theorem~\ref{tda1} and we have $z_1=\tilde{y}$. Now we would like $z_2$ to be $\tilde{y}'$.
So we differentiate $\frac{1}{\sqrt{r}}(-\theta_0 y + y')$ modulo the relation ${\mathcal{L}}_my=0$: this gives us $\tilde{y}'$ as a linear combination of $y$ and $y'$ and we find the expression of $P_m$
giving $\tilde{X}=P_m X$.

The new coefficient field is $\widetilde{K}:=K(\theta_0)$.
When $\theta_0$ is an algebraic function over $K$,
$\widetilde{K}$ is an algebraic extension of the differential field $K$ and the Picard--Vessiot extension $\widetilde{F}_m$ of equation $\widetilde{\mathcal L}_m u=0$ is an algebraic extension of the Picard--Vessiot extension $F_m$ of ${\mathcal{L}}_my=0$.
Thus, $(\mathrm{DGal}(F_m/K))^\circ =\big(\mathrm{DGal}\big(\widetilde F_m/K\big)\big)^\circ$.
Finally, if $\theta_0$ belongs to $K$, then $\widetilde K=K$ and $\widetilde{F}_m=F_m$, which implies that the Darboux transformation preserves the Galois groups.
\end{proof}

\begin{Remark}
This matrix factorization has an interesting interpretation. The gauge transformation is $\widetilde{X} = P_m X$.
Now, it is easily seen that
\begin{gather*}
R X = \begin{pmatrix} \frac{1}{\sqrt{r}}y \\ \tilde{y} \end{pmatrix}\!.
\end{gather*}
Indeed, see Theorem~\ref{tda1}, the second row is the Darboux transformation. The matrix factorization $P_m=L_m R$, combined with this relation, thus shows
that
\begin{gather*}
\begin{pmatrix} \tilde{y} \\ \tilde{y}' \end{pmatrix} = L_m \begin{pmatrix} \frac{1}{\sqrt{r}}y \\ \tilde{y} \end{pmatrix}\!.
\end{gather*}
The first row is an obvious identity.
The second row provides us a notably simple first-order link between solutions $y$ of $\mathcal{L}_m$ and solutions $\tilde{y}$ of its Darboux transformation $\widetilde{\mathcal{L}}_m$:
\begin{gather*}
\tilde{y}' - \rho \tilde{y} = m \sqrt{r} y.
\end{gather*}
Note that $\rho$ does not depend on the parameter: $m$ appears only in the right hand side.
This observation provides us with a simple explicit form for the reverse of the Darboux transformation. Given a solution $\tilde{y}$ of the transformed equation, we may transform it back into a solution $y$ of the initial equation via $y = \frac1{m \sqrt{r}}(\tilde{y}' - \rho \tilde{y})$.
\end{Remark}

\subsection{Systems in \texorpdfstring{\mathversion{bold}$\mathfrak{s}\mathrm{y}\mathfrak{m}^2(\mathrm{SL}(2,C_K))$}{sym\string^2(SL(2,C\_K))} and \texorpdfstring{$\mathfrak{so}(3,C_K)$}{so(3,C\_K)}}

Now we use the method of Section~\ref{tensor-constructions} to build linear differential systems in $\sym^2 (\mathrm{SL}(2,C_K))$ and give their relations to systems in $\mathfrak{so}(3,C_K)$.

Consider the linear differential system $[A_0]$ as in equation~\eqref{mdar1} for $m=0$, where
$p=w'/w$ and $w,q\in K$. We recall that its second symmetric power system is given by the linear differential system
\begin{gather}\label{ssyp}
[\sym^2 (A_0)]\colon\ Y'=-S_2 Y,
\end{gather}
for
\begin{gather*}
Y_2:= \Sym^2 (X) = \begin{pmatrix} y^2\\2yy'\\(y')^2 \end{pmatrix}
\quad \textrm{ and } \quad S_2:= \sym^2 (A_0) = \begin{pmatrix} 0&-1&0\\2q &p
	&-2\\0&q& 2p
	\end{pmatrix}\!.
\end{gather*}
We recall that we showed, in Example~\ref{ex-second-sym}, how to build a solution matrix and a representation of the Galois group in this construction.
We record that for further use in this classical lemma.

\begin{Lemma}\label{lem:gal group}
Let $\mathbf{X}$ be a fundamental matrix of system~\eqref{mdar1}
and $\mathrm{DGal}(F_m/K)$
be
its dif\-fer\-en\-tial Galois group. Then, $\mathbf{Y}= \Sym^2 (\mathbf{X})$ is a fundamental matrix for system~\eqref{ssyp} and $\Sym^2 (\mathrm{DGal}(F_m/K))$ is a representation of its differential Galois group.
\end{Lemma}

This situation can be illustrated by the following diagram:
\begin{gather}\label{eq:diag 1}
\begin{split}
 \xymatrix{ [A]\colon X'= -A X \ar@{~>}[rr] \ar[d]^{\sym^2} & & \mathrm{DGal}(F_0/K) \ar[d]^{\Sym^2} \\ [ \sym^2 (A)]\colon Y'=-\sym^2 (A) Y \ar@{~>}[rr] & & \Sym^2 (\mathrm{DGal}(F_0/K)). }
 \end{split}
\end{gather}

The next result provides us with a gauge transformation to go from an $\mathfrak{so}(3, C_K)$ system
to a~second symmetric power system of the form~\eqref{ssyp}. This enables us to express equation $\mathcal L_0y=0$ as a linear differential system in $\mathfrak{so}(3, C_K)$ written in the form $Z'= Z \times \Omega$ (equation~\eqref{fedo}). As~a~consequence of this result, we can extend previous reasoning for differential Galois groups to $\mathfrak{so}(3, C_K)$ systems, as we will show next.

\begin{Lemma}\label{gaugemat}
Let $Q$ be the gauge matrix given by
\begin{gather}\label{eq:gauge sltoso3}
 Q=
 \begin{pmatrix}
 	1&\hphantom{-}0&-1\\ {\rm i}&\hphantom{-}0&\hphantom{-}{\rm i}
\\ 0&-1&\hphantom{-}0
	\end{pmatrix}\!.
\end{gather}
The gauge transformation $Z=QY$ transforms the $\mathfrak{so}(3, C_K)$ system $Z'= Z \times \Omega$, where $\Omega=(f,g,h)^{\rm T}$ into the system $Y'=\sym^2(C)Y$
with
\begin{gather*}
C:= \frac{1}{2} \begin{pmatrix} {\rm i} h & (g +{\rm i}f ) \\
	 - (g - {\rm i}f ) & - {\rm i} h \end{pmatrix}\!.
\end{gather*}
\end{Lemma}

\begin{proof}
We have
\begin{gather*} \sym^2(C) = \begin{pmatrix}
{\rm i}h &\frac12(g +{\rm i} f) &0\\ -(g-{\rm i}f) &0&g+{\rm i}f \\ 0&-\frac12(g -{\rm i} f) &-{\rm i}h \end{pmatrix}
\end{gather*}
As $Q$ is constant, the effect of the gauge transformation is just conjugation by $Q$.
We have
\begin{gather*}
Q.\sym^2(C) = \begin{pmatrix}
{\rm i}h &g &{\rm i}
h \\ -h &-f
 &h \\ g -{\rm i}f &0&-g -{\rm i}f \end{pmatrix}
\quad \textrm{and} \quad
Q.\sym^2(C).Q^{-1} =\begin{pmatrix} 0&h &-g
\\ -h &0&f
\\ g &-f &0
\end{pmatrix}\!.\!\!\!\tag*{\qed}
\end{gather*}
\renewcommand{\qed}{}
\end{proof}

This shows that, given an orthogonal system, we have an explicit gauge transformation formula to view it as a symmetric square of a second-order system.

\begin{Remark}
Using an additional gauge transformation, one can transform the matrix $C$ of the above Lemma~\ref{gaugemat}
into a companion form. This allows one to recover the formulas from Remark~\ref{rmk-darboux-thm2} and reprove Theorem~\ref{tso3dar}.

Conversely, if we start from $ \mathcal{L}_0 y:= y''+ p y'+ q y=0$, reversing the above process will produce an equivalent orthogonal system $Z'=Z \times \Omega$ with
\begin{gather*}
\Omega=( {\rm i} (q-1) , q+1,{\rm i} p)^{\rm T} =: (f,g,h)^{\rm T}.
\end{gather*}
\end{Remark}

In~\cite{Blaz-M}, see also~\cite{Blaz}, D.~Bl\'azquez-Sanz and J.J.~Morales-Ruiz have also given such a classical Lie algebra isomorphism between $\mathfrak{so}(3, C_K)$ and $\mathfrak{sl} (2, C_K)$, see~\cite[Proposition~6.7]{Blaz-M} (a small matrix calculation shows that they satisfy the same bracket rules):
\begin{gather*}
 \begin{pmatrix} & 1 & \\ -1 & & \\ & & 0 \end{pmatrix} \mapsto \begin{pmatrix} \frac{\rm i}{2} & 0 \\ 0 & -\frac{\rm i}{2} \end{pmatrix}\!, \qquad
 \begin{pmatrix} & & 1 \\ & 0 & \\ -1 & & \end{pmatrix} \mapsto \begin{pmatrix} 0 & \frac{1}{2} \\ -\frac{1}{2} & 0 \end{pmatrix}\!, \\
 \begin{pmatrix} 0 & & \\ & & 1 \\ & -1 & \end{pmatrix} \mapsto \begin{pmatrix} 0 & -\frac{\rm i}{2} \\ -\frac{\rm i}{2} & 0 \end{pmatrix}\!.
\end{gather*}

\begin{Corollary}\label{corso3}
A fundamental matrix for system $Z' = Z \times \Omega$ is
\begin{gather}\label{eq:mat ftal Z}
 \mathbf{Z}:= Q \mathbf{Y} = w \begin{pmatrix}
y_1^2 - (y_1')^2 & y_1 y_2 - y_1' y_2' & y_2^2 -(y_2')^{2} \\[3pt]
{\rm i} \big(y_1^2 + (y_1')^{2}\big) & {\rm i} (y_1 y_2 + y_1' y_2') & {\rm i} \big(y_2^2 + (y_2')^{2}\big) \\[3pt]
-2 y_1 y_1' & - y_1 y_2' - y_1' y_2 & -2 y_2 y_2' \end{pmatrix}
\end{gather}
for $Q$ and $\mathbf{Y}$ defined by~\eqref{eq:gauge sltoso3} and~\eqref{eq:mat ftal Y}, respectively.
\end{Corollary}

Finally, we can compute the differential Galois group of system~\eqref{fedo}.

\begin{Corollary}
Using the fundamental matrix~\eqref{eq:mat ftal Z}, the matrices in the differential Galois group of the $\mathfrak{so} (3)$ system $Z' = Z \times \Omega$ are the matrices $\Sym^2 (M_{\sigma})$ of Lemma~$\ref{lem:gal group}$.
\end{Corollary}

\begin{proof}
Let $\sigma$ be in the Galois group. It acts on $\mathbf{Y}$ via $\sigma (\mathbf{Y}) = \mathbf{Y} \cdot \Sym^2 (M_{\sigma})$ (Lemma~\ref{lem:gal group}). Now, $\sigma (Q)= Q$ because $w \in K$. So, we have $\sigma (\mathbf{Z}) = \sigma (Q \mathbf{Y}) = Q \cdot \mathbf{Y} \cdot \Sym^2 (M_{\sigma}) = \mathbf{Z} \cdot \Sym^2 (M_{\sigma})$.
\end{proof}

Now, we can complete diagram~\eqref{eq:diag 1} by adding the action on the $\mathfrak{so}(3, C_K)$ system:
\begin{gather}\label{eq:diag 2}
\begin{split}
\xymatrix{ [A_0]\colon X'= -A_0 X \ar@{~>}[rr] \ar@<1ex>[d]^{\sym^2} & & \mathrm{DGal}(L_0/K) \ar@<1ex>[d]^{\Sym^2}
\\
[ \sym^2 (A_0)]\colon Y'=-S_2 Y \ar@{~>}[rr] \ar@<1ex>[d]^{Q} \ar[u]
 & & \Sym^2 (\mathrm{DGal}(L_0/K))
\ar@{=}[d]
\\
[\Omega]\colon Z' = - \Omega Z \ar@{~>}[rr] \ar[u] & & \Sym^2 (\mathrm{DGal}(L_0/K)). }
\end{split}
\end{gather}

As noted by a referee, the arrow from $\Sym^2 (\mathrm{DGal}(L_0/K)) $ back to $\mathrm{DGal}(L_0/K)$ is not well defined because, for any matrix $M$, $\Sym^2(M) = \Sym^2(-M)$. However, the arrows between systems (on the left) are well defined and this is all that we need for our constructions.
Finally, we obtain a simple corollary.
\begin{Corollary}
Consider the orthogonal system $[\Omega]\colon Z' = \Omega \times Z$ with $Z = (\alpha , \beta , \gamma)^{\rm T}$ and the equivalent second symmetric power $\big[\sym^2 (A_0)\big]\colon Y' = \sym^2 (A_0) Y$ with $Y= (z_1, z_2 , z_3)^{\rm T}$. The system $[\Omega]$ admits the first integral $\alpha^2 + \beta^2 + \gamma^2$ and $\big[\sym^2 (A_0)\big]$ admits the first integral $w^2 \big(4 z_1 z_3 - z_2^2\big)$.
\end{Corollary}

\begin{proof}The first part is well known and is proved by Darboux in~\cite[Part~I, Chapter~I, p.~8 and Chapter~II, p.~28]{da3}.
The second part follows from the application of the gauge transformation of Lemma~\ref{gaugemat}: it transforms the first integral $\alpha^2 + \beta^2 + \gamma^2$ of $[\Omega]$ into $w^2 \big(4 z_1 z_3 - z_2^2\big)$:
this is still a constant of motion and hence a first integral of $\big[\sym^2 (A_0)\big]$.
\end{proof}

\subsection{Darboux transformations in \mathversion{bold}\texorpdfstring{$\Sym^2(\mathrm{SL}(2,C_K))$}{Sym\string^2(SL(2,C\_K))} and \texorpdfstring{$\mathfrak{so}(3,C_K)$}{so(3,C\_K)}}

The aim of this subsection is to construct Darboux transformations for third-order orthogonal systems using diagram~\eqref{eq:diag 2}.
Our construction will ensure that these Darboux transformations will preserve the identity component of the differential Galois group of each equation.

In order to construct the Darboux transformations for the second symmetric power system coming from the general second-order linear differential equation ${\mathcal{L}}_my=0$ we extend Proposition~\ref{prop1}. We present two ways to extend it. Like in previous subsection, this will allow us to obtain Darboux transformations for the $\mathfrak{so}(3, C_K)$ system.

For the first one, consider the linear differential system $[A_m]\colon X'=-(A_0+mN) X$ given by equation~\eqref{mdar1}. Its second symmetric power system is given by the linear differential system
\begin{gather}\label{ssyp m}
\big[\sym^2 (A_m)\big]\colon\ Y'=-(S_2 +m N_2) Y ,
\end{gather}
where $Y=\Sym^2 (X) = \big(y^2 ,2yy', (y')^2\big)^{\rm T}$ and $S_2 +m N_2 = \sym^2 (A_0 +mN)$ are given by
\begin{gather*}
S_2:=\begin{pmatrix} 0&-1&0\\2q &p &-2\\0&q&2p\end{pmatrix} \qquad \text{and} \qquad N_2:=\begin{pmatrix}0&0&0\\ -2r& 0& 0\\0&- r&0\end{pmatrix}\!.
\end{gather*}
A fundamental matrix for this system is given by matrix~\eqref{eq:mat ftal Y}, where $\{ y_1, y_2 \}$ is a basis of solutions of equation~\eqref{eqdarb}.

Recall that, after applying the gauge transformation~\eqref{gauge}, system~\eqref{mdar1} is transformed into the linear differential system $\big[\widetilde{A}_m\big]\colon \widetilde{X}'=-\big(\widetilde A_0+mN\big) \widetilde{X}$, defined by equation~\eqref{mdar2}, whose second symmetric power system is given by the linear differential system
\begin{gather}\label{ssyp m dt}
\big[\sym^2 (\widetilde{A}_m)\big]\colon\ \widetilde{Y}'=-\big(\widetilde{S}_2+ m N_2\big) \widetilde{Y} ,
\end{gather}
where $\widetilde{Y}= \Sym^2 \big(\widetilde{X}\big) = \big(\tilde{y}^2 ,2\tilde{y} \tilde{y}', (\tilde{y}')^2\big)^{\rm T}$ and
$\widetilde{S}_2+ m N_2 = \sym^2 \big(\widetilde{A}_0 +mN\big)$, for
\begin{gather*}
\widetilde{S}_2:=\begin{pmatrix} 0&-1&0 \\ 2\tilde{q} &p &-2 \\ 0& \tilde{q} &2p\end{pmatrix}\!,
\end{gather*}
where $\tilde{y}$ and $\tilde{q}$ have the explicit form expressed in Theorem~\ref{tda1}. Thus, since $\widetilde{Y}= \Sym^2 \big(\widetilde{X}\big) = \Sym^2 (P_m) \cdot \Sym^2 (X) = \Sym^2 (P_m) \cdot Y $, the gauge transformation~\eqref{gauge} also induces a transformation in the second symmetric power systems which sends system~\eqref{ssyp m} to system~\eqref{ssyp m dt}.
The following result formalizes this idea.

\begin{Proposition}[first Darboux transformation for $\Sym^2 (\mathrm{SL}(2,C_K))$] \label{thnudarsl21}
Let $P_{1,m}:= \Sym^2 (P_{m})$.
We have 
\begin{gather} \label{eq:DT sl2}
P_{1,m} = \dfrac{1}{r}
\begin{pmatrix} {\theta_0}^{2}&-\theta_0&1 \\ 2 \theta_0
\nu & \nu-\theta_0 \rho& 2 \rho \\
\nu^{2}& -\rho \nu&{\rho}^{2}
 	\end{pmatrix} 	=
	 \frac{1}{r} \begin{pmatrix} 0&0&1\\ 0&mr&2 \rho
\\ {m}^{2}{r}^{2}&\rho mr&{\rho}^{2}\end{pmatrix} \cdot
\begin{pmatrix} 1&0&0\\ -2 \theta_{{0}}
&1&0\\ \theta_{{0}}^{2}&-\theta_{{0}}&1
\end{pmatrix}\!,
\end{gather}
where $P_m$ is defined by expression~\eqref{gauge},
	$\theta_0= \frac{y_0'}{y_0}$, $\rho = -\theta_0- p - \frac{r'}{2r}$
	and $\nu=mr-\theta_0 \rho$.
Then, $P_{1,m}$ is a gauge transformation which sends system $\sym^2 ([A_m])$ to system $\sym^2 \big(\big[\widetilde{A}_m\big]\big)$.
\end{Proposition}

\begin{proof}Since matrix $P_{1,m}$ is the second symmetric power matrix of matrix $P_m$, defined by~\eqref{gauge}, given $Y$ and $\widetilde{Y}$ solutions of equations~\eqref{ssyp m} and~\eqref{ssyp m dt} respectively, it satisfies $\widetilde{Y} = \Sym^2 (P_{m}) Y$.
 Proposition~\ref{prop1} shows that $P_{m} = L_m. R$. As $\Sym(\bullet)$ is a group morphism, we have
$\Sym^2(P_m)=\Sym^2(L_m)\Sym^2(R)$, which gives the matrix factorization and the result.
\end{proof}

As systems $\sym^2 ([A_m])$ and $\sym^2 \big(\big[\widetilde{A}_m\big]\big)$ have the same shape, we may say that $P_{1,m}$ is a~\emph{Dar\-boux transformation} from $\sym^2 ([A_m])$ to $\sym^2 \big(\big[\widetilde{A}_m\big]\big)$.
This new Darboux transformation is induced by the original Darboux transformation for second-order systems. We recover the matrix factorization from our expression of the Darboux transformation as a gauge transformation in Proposition~\ref{prop1}.

There is another way to build a Darboux transformation for second symmetric power systems. That is, transform the first system $[A_m]\colon X'=-(A_0+mN) X$, given by equation~\eqref{mdar1}, into a system in $\mathfrak{sl} (2,K)$ and perform all the previous process with the resulting system. By doing that, we ensure that the differential Galois group of equation~\eqref{mdar1}
is in ${\rm SL}(2, C_K)$.
In order to obtain an $\mathfrak{sl} (2,K)$ system, we consider the gauge change:
\begin{gather} \label{eq:gauge sl2}
X_{1}:=\Delta X , \qquad
\Delta:= \begin{pmatrix} 1&0 \\ 0 & w \end{pmatrix}\!,
\end{gather}
where $w$ denotes the inverse of the Wronskian of $[A_m]$
(it is a non-zero solution of $w'=\textrm{trace}(A_0+mN).w$). We have
$X_{1} = (y, {w}y')^{\rm T} $ and the resulting system is the $\mathfrak{sl} (2,K)$ system
\begin{gather}\label{eq:sist sl2}
 [B_{m}]\colon X_{1}'=-(B_0+mN_1) X_{1} ,
\end{gather}
with $B_0 + mN_1 \in \mathfrak{sl}(2,K)$, given by
\begin{gather*}
B_0:= \begin{pmatrix} 0& - \frac{1}{w}\\ wq&0 \end{pmatrix} \qquad \text{and} \qquad
N_1:= \begin{pmatrix} 0&0\\ -wr&0 \end{pmatrix}\!.
\end{gather*}
By an elimination process (cyclic vector, see end of Section~\ref{dgal-section}), we can show that this systems is still equivalent to the second-order linear differential equation~\eqref{eqdarb}.
And a fundamental matrix for this system is given by
\begin{gather}\label{eq:sol-fund-gaug-2}
 \mathbf{X}_{1} = \begin{pmatrix} y_1 & y_2 \\ wy_1' & w y_2' \end{pmatrix}\!,
\end{gather}
where $\{ y_1, y_2 \}$ is a basis of solutions of equation~\eqref{eqdarb}.

The Darboux transformation $P_m$, given by~\eqref{gauge}, for system~\eqref{mdar1} induces the Darboux transformation $\Delta P_m \Delta^{-1}$ for system~\eqref{eq:sist sl2}.

Now, we consider the second symmetric power system of system~\eqref{eq:sist sl2}. This system is given~by
\begin{gather}\label{eq:sist ssmsl2}
\big[\sym^2 (B_m)\big]\colon \ Y_1'=-\big(\widehat{S}_2 + m \widehat{N}_2\big) Y_1 ,
\end{gather}
where $Y_1:= \Sym^2 ( X_1) = \big( y_1^2, 2y_1y_2 , y_2^2 \big)^{\rm T}$ and $\widehat{S}_2 + m \widehat{N}_2 =\sym^2 (B_0 +mN_1) \in \mathfrak{sl} (3,K)$, for
\begin{gather*}
\widehat{S}_2:=\begin{pmatrix} 	
0&- \frac{1}{w} &0\\ 2 wq&0&
\frac{-2}{w}\\ 0&wq&0
	\end{pmatrix} \qquad \text{and} \qquad
	\widehat{N}_2:=\begin{pmatrix}
	0&0&0\\ -2 wr&0&0
\\ 0&-wr&0
	\end{pmatrix}\!.
\end{gather*}
From the above, we can easily compute the expression of a fundamental matrix for system~\eqref{eq:sist ssmsl2}:
\begin{gather}\label{eq:sol-fund-sl-2}
 \mathbf{Y}_{1} = \Sym^2 ( \mathbf{X}_1) = \begin{pmatrix} y^2_1 & y_1 y_2 & y_2 \\[2pt] 2wy_1 y_1' & w(y_1 y_2' + y_1' y_2) & 2wy_2 y_2' \\[3pt] w^2 (y_1')^2 & w^2 y_1' y_2' & w^2 (y_2')^2 \end{pmatrix}\!,
\end{gather}
where $\mathbf{X}_1$ is given by~\eqref{eq:sol-fund-gaug-2}.

Next, we find the second expression for the Darboux transformation for systems that can be written as a second symmetric power.

\begin{Corollary}[second Darboux transformation for $\Sym^2 (\mathrm{SL}(2,C_K))$] \label{thnudarsl22} Consider the system $\big[\sym^2 (B_m)\big]\colon Y_1'=-\big(\widehat{S}_2 + m \widehat{N}_2\big) Y_1$ given by~\eqref{eq:sist ssmsl2}.
Let $P_{2,m} :=\Sym^2 \big(\Delta P_m \Delta^{-1}\big)$. We have
\begin{align}
 P_{2,m} &= \dfrac{1}{r} \begin{pmatrix} {\theta_0}^{2}&-w \theta_0&{w}^{2}
\\ \frac{-2 \theta_0 \nu }{w}&\nu- \rho \theta_0&2 w \rho
\\ \frac {\nu^{2}}{{w}^{2}}&{\frac {\rho \nu }{w}}&{\rho}^{2}\end{pmatrix}\nonumber
\\
&= \frac1{r}\begin{pmatrix}
0&0&1\\ 0&wmr&2 w\rho
\\ {w}^{2}{m}^{2}{r}^{2}&{w}^{2}\rho mr&{w}^{2}{
\rho}^{2}
\end{pmatrix}
\cdot
\begin{pmatrix}
1&0&0\\ -2 \theta_{{0}}
&\frac1{w}&0\\ \theta_{{0}}^{2}&-{\frac{\theta_{0}}{w}}& 	\frac1{w^2}
\end{pmatrix}\!,\label{eq:DT sl22}
\end{align}
where matrices $P_m$ and $\Delta$ are defined by~\eqref{gauge} and~\eqref{eq:gauge sl2}, respectively, $\theta_0= \frac{y_0'}{y_0}$, $\rho = -\theta_0- p - \frac{r'}{2r}$
and $\nu=mr - \rho \theta_0$.
Then, $P_{2,m}$ is a Darboux transformation for system $\big[\sym^2 (B_m)\big]$.
\end{Corollary}

\begin{proof}
As we have seen, a Darboux transformation for the $\mathfrak{sl} (2,K)$ system~\eqref{eq:sist sl2} is given by $\widetilde{X}_1 = \Delta P_m \Delta^{-1} X_1$, where $\Delta$ and $P_m$ are defined by~\eqref{eq:gauge sl2} and~\eqref{gauge}, respectively. The transformed system is the $\mathfrak{sl} (2,K)$ system
\begin{gather*}
 \big[\widetilde{B}_{m}\big]:= \widetilde{X}_{1}'=-\big(\widetilde{B}_0+mN_1\big) \widetilde{X}_{1},
\end{gather*}
with $\widetilde{B}_0+mN_1 \in \mathfrak{sl} (2,K)$, for
\begin{gather*}
\widetilde{X}_{1} = \begin{pmatrix} \tilde{y}_1 \\ \tilde{y}_2 \end{pmatrix} \qquad \text{and} \qquad
\widetilde{B}_0:= \begin{pmatrix} 0 & \frac1w \\ -\tilde{q} w & 0 \end{pmatrix}\!.
\end{gather*}
Its corresponding second symmetric power system is
\begin{gather}\label{ssyp sl2}
\big[\sym^2 \big(\widetilde{B}_m\big)\big]\colon\ \widetilde{Y}_1'=-\big(\widetilde{\widehat{S}}_2 + m \widehat{N}_2\big) \widetilde{Y}_1,
\end{gather}
where $ \widetilde{Y}_1:= \Sym^2 \big(\widetilde{X}_1\big) = \big( \tilde{y}_1^2, 2\tilde{y}_1\tilde{y}_2 , \tilde{y}_2^2 \big)^{\rm T}$
and $\widetilde{\widehat{S}}_2 + m \widehat{N}_2 =\sym^2 \big(\widetilde{B}_0+mN_1\big) \in \mathfrak{sl} (3,K)$, for
\begin{gather*}
\widetilde{\widehat{S}}_2:=\begin{pmatrix}
0&- \frac{1}{w}&0 \\ 2 w \tilde{q} &0& \frac{-2}{w} \\ 0&w \tilde{q}&0
\end{pmatrix}\!,
\end{gather*}
where $\tilde{q}$ is the ``new potential'' obtained by the Darboux transformation in Theorem~\ref{tda1}.

Now, consider the second symmetric power system~\eqref{eq:sist ssmsl2}. Since
\begin{gather*}
\widetilde{Y}_1 = \Sym^2 \big(\widetilde{X}_1\big) = \Sym^2 \big(\Delta P_m \Delta^{-1}\big) \cdot \Sym^2 ( X_1) = \Sym^2 \big(\Delta P_m \Delta^{-1}\big) \cdot Y_1,
\end{gather*}
we find that matrix $P_{2,m}:=\Sym^2 \big(\Delta P_m \Delta^{-1}\big)$
is a Darboux transformation which sends system~\eqref{eq:sist ssmsl2} into system~\eqref{ssyp sl2}.
We had the matrix factorization $P_m=L_m R$ (Proposition~\ref{prop1}). So,
as $\Sym(\bullet)$ is a group morphism,
\begin{gather*}
\Sym^2 \big(\Delta P_m \Delta^{-1} \big) = \Sym^2(\Delta L_m)\cdot \Sym^2\big(R \Delta^{-1}\big)
 \end{gather*}
and this gives us the desired matrix factorization above.
\end{proof}

We recall that if $w=1$, Darboux transformations $P_{1,m}$ and $P_{2,m}$ are the same. This is useful in the applications to non-relativistic and one-dimensional quantum mechanics where $r=1$ as well, see Section~\ref{sec:SUSY}.

Once we have defined the Darboux transformations for second symmetric power systems, we can state the Darboux transformation for $\mathfrak{so}(3, C_K)$ systems as follows. As we have found two Darboux transformations for second symmetric power systems, we will have two Darboux transformations for $\mathfrak{so}(3, C_K)$ systems as well: one using Proposition~\ref{thnudarsl21} and another one using
Corollary~\ref{thnudarsl22}.

Applying Lemma~\ref{gaugemat}, we can transform the second symmetric power system~\eqref{ssyp m} into the $\mathfrak{so}(3, C_K)$ system
\begin{gather}\label{eq:system so3}
 [\Omega_m]\colon\ Z' = - (\Omega_0 + mN_3 ) Z,
\end{gather}
where $Z= Q Y$ and $- (\Omega_0 + mN_3 ) = Q' Q^{-1} - Q (S_2 + m N_2) Q^{-1} $ are given by
\begin{gather*}
Z= \begin{pmatrix}
	\alpha \\ \beta \\ \gamma \end{pmatrix}\!, \qquad
	\Omega_0=\begin{pmatrix}
	0&-{\rm i}p&-(q+1)\\ {\rm i}p&0&{\rm i}(q-1) \\ q+1&-{\rm i}(q-1)&0
	\end{pmatrix} \qquad \text{and} \qquad
	N_3 = r \begin{pmatrix} 0 & 0 & -1 \\ 0& 0& {\rm i} \\ 1 & -{\rm i} & 0 \end{pmatrix}\!.
\end{gather*}
A fundamental matrix for this system is given by matrix~\eqref{eq:mat ftal Z}, where $\{y_1, y_2\}$ is a basis of solutions of equation~\eqref{eqdarb}.

After performing the Darboux transformation~\eqref{eq:DT sl2}, system~\eqref{ssyp m} is transformed into system~\eqref{ssyp m dt}, which, again by Lemma~\ref{gaugemat}, can be transformed into the $\mathfrak{so}(3, C_K)$ system
\begin{gather}\label{eq:system so3 DT}
 \big[\widetilde{\Omega}_m\big]\colon\ \widetilde{Z}' = - \big(\widetilde{\Omega}_0 + m N_3 \big) \widetilde{Z},
\end{gather}
where $\widetilde{Z}= Q \widetilde{Y}$ and $- \big(\widetilde{\Omega}_0 + m N_3 \big) = Q' Q^{-1} - Q \big(\widetilde{S}_2 + m N_2\big) Q^{-1} $ for
\begin{gather*}
\widetilde{Z}= \begin{pmatrix} \widetilde{\alpha} \\ \widetilde{\beta} \\ \widetilde{\gamma} \end{pmatrix}
\qquad \text{and} \qquad
\widetilde{\Omega}_0=
\begin{pmatrix}
0&-{\rm i}p&-(\tilde{q}+1)\\ {\rm i}p&0&{\rm i}(\tilde{q}-1) \\ \tilde{q}+1&-{\rm i}(\tilde{q}-1)&0
\end{pmatrix}\!,
\end{gather*}
for $\tilde{q}$ as in Theorem~\ref{tda1}.
Thus, the Darboux transformation~\eqref{eq:DT sl2} also induces a transformation in the corresponding $\mathfrak{so}(3, C_K)$ systems which sends system~\eqref{eq:system so3} to system~\eqref{eq:system so3 DT}.

The following result shows that this induced transformation is indeed a Darboux transformation for $\mathfrak{so}(3, C_K)$ systems.

\begin{Proposition}[first Darboux transformation for $\mathfrak{so}(3, C_K)$]\label{secthdars031}
Consider the systems $[\Omega_m]$\emph{:} $Z' = - (\Omega_0 + mN_3 ) Z$ and $\big[\widetilde{\Omega}_m\big]\colon \widetilde{Z}' = - \big(\widetilde{\Omega}_0 + m N_3 \big) \widetilde{Z}$ given by~\eqref{eq:system so3} and~\eqref{eq:system so3 DT}, respectively.
 Let $T_{1,m}$ be the matrix defined by
\begin{align*}
 T_{1,m} = Q P_{1,m} Q^{-1} =\frac{1}{2 r} \begin{pmatrix}
-{\nu}^{2}+{\rho}^{2}+{\theta_0 }^{2}-1 & {\rm i} ({\nu}^{2}+{\rho}^{2}-{\theta_0}^{2}-1) & 2 (\nu \rho+\theta_0)
\\
{\rm i} \big({\nu}^{2}-{\rho}^{2}+{\theta_0 }^{2}-1\big) &{\nu}^{2}+{\rho}^{2}+{\theta_0 }^{2}+1 & 2 {\rm i} (\theta_0 -\nu \rho)
\\
2 ( \nu \theta_0 +\rho) &-2 {\rm i} (\nu \theta_0 -\rho) &2 ( \nu - \theta_0 \rho) \end{pmatrix}\!,
\end{align*}
where matrix $Q $ is defined by~\eqref{eq:gauge sltoso3},
matrix $P_{1,m}$ is defined by expression~\eqref{eq:DT sl2} and $\nu=mr-\theta_0 \rho$.
Then, $ T_{1,m}$ is a Darboux transformation, i.e., a
gauge transformation which sends system $ [\Omega_m]\colon Z' = - (\Omega_0 + mN_3 ) Z$ to a system $\big[\widetilde{\Omega}_m\big]\colon \widetilde{Z}' = - \big(\widetilde{\Omega}_0 + m {N}_3 \big) \widetilde{Z}$
of the same shape: we still have $\widetilde{\Omega}_0\in \mathfrak{so}(3)$, the spectral part $N_3$ is untouched and the Galoisian properties of the system are preserved.
\end{Proposition}

\begin{Remark}
As in the previous results, the matrix can be factored into an $m$-dependent part and an independent part:
\begin{gather*}
T_{1,m} = \begin{pmatrix} -{m}^{2}{r}^{2}&-\rho mr&-{\rho}^{2}+1
\\ {\rm i}{m}^{2}{r}^{2}&{\rm i}\rho mr&{\rm i}+{\rm i}{\rho}^{2}
\\ 0&-mr&-2 \rho
\end{pmatrix}\!. \begin{pmatrix} \frac{1}{2} &-\frac{\rm i}{2}&0\\ -\theta_{{0}} &{\rm i} \theta_{{0}}&-1\\
\frac{1}{2}\big(\theta_0^{2} -1\big)
&
-\frac{\rm i}{2} \big(\theta_0^{2}+1\big)
&\theta_{{0}}
\end{pmatrix}\!,
\end{gather*}
where the $m$-dependent matrix is
$Q\, \Sym^2(\Delta L_m)$ and the matrix on the right is $\Sym^2\!\big(R \Delta^{-1}\big) Q^{-1}$.
\end{Remark}

\begin{proof}
The proof follows from the application of Lemma~\ref{gaugemat} and Theorem~\ref{thnudarsl21}. Given $Y$ and $\widetilde{Y}$ solutions of equations~\eqref{ssyp m} and~\eqref{ssyp m dt} respectively, and $Z$ and $\widetilde{Z}$ solutions of~\eqref{eq:system so3} and~\eqref{eq:system so3 DT}, respectively, by Lemma~\ref{gaugemat}, we have that $Z= Q Y$ and $\widetilde{Z} = Q \widetilde{Y}$. On the other hand, by Theorem~\ref{thnudarsl21}, we know that $\widetilde{Y} = P_{1,m} Y$. Thus, we get the gauge transformation $\widetilde{Z} = \big(Q P_{1,m} Q^{-1}\big) Z = T_{1,m} Z$.
From this, we immediately obtain the gauge transformation for the coefficient matrix:
\begin{gather*}
- \big(\widetilde{\Omega}_0 + m N_3 \big) = T'_{1,m} T^{-1}_{1,m} - T_{1,m}(\Omega_0 + mN_3) T^{-1}_{1,m} .
\end{gather*}

The rest of the corollary is proved following the same argument as in Proposition~\ref{prop1}.
\end{proof}

Propositions~\ref{thnudarsl21} and~\ref{secthdars031} can be summarized in the following commutative diagram:
\begin{gather}\label{eq:diag-1}
\begin{split}
\xymatrix{ [A_m]\colon X'= -(A_0 +mN) X \ar[rr]^{P_m} \ar@<1ex>[d]^{\sym^2} & & [\widetilde{A}_m]\colon \widetilde{X}'= -\big(\widetilde{A}_0 +mN\big) \widetilde{X} \ar@<1ex>[d]^{\sym^2} \\ \big[\sym^2 (A_m)\big]\!:= Y'\!={-}(S_2\! +\!m N_2) Y\! \ar[rr]^{P_{1,m}} \ar@<1ex>[d]^{Q } \ar[u] & & \!\big[\sym^2 \big(\widetilde{A}_m\big)\big]\! :=\widetilde{Y}'\!={-}\big(\widetilde{S}_2\! +\!m N_2\big) \widetilde{Y} \ar@<1ex>[d]^{Q} \ar[u] \\ [\Omega_m]\colon Z' = - (\Omega_0 + mN_3 ) Z \ar[rr]^{ T_{1,m}} \ar[u] & & [\widetilde{\Omega}_m]\colon \widetilde{Z}' = - \big(\widetilde{\Omega}_0 + m N_3 \big) \widetilde{Z} \ar[u]}.
\end{split}\hspace{-10mm}
\end{gather}

We end this section by establishing a second Darboux transformation for $\mathfrak{so}(3, C_K)$ systems. For that, we transform the $\Sym^2 ({\rm SL}(2,K))$ system~\eqref{eq:sist ssmsl2} into an
${\rm SO}(3, C_K)$ system. Consider the matrix
\begin{gather*}
S = \begin{pmatrix} 1&0&\hphantom{-}1 \\ 0&{\rm i}&\hphantom{-}0\\ {\rm i}&0&-{\rm i}\end{pmatrix}
\end{gather*}
and the gauge change $Z_1= S \cdot Y_1 $. Then, the system
\begin{gather}\label{eq:so-system-2}
 \big[\widehat{\Omega}_m\big]\colon \ Z_1' = - \big(\widehat{\Omega}_0 + m \widehat{N}_3 \big) Z_1,
\end{gather}
where
\begin{gather*}
 \widehat{\Omega}_0 = \begin{pmatrix}
 0 & {\rm i} \big(\frac{1}{w} - wq \big) & 0 \vspace{1mm}\\
 -{\rm i} \big(\frac{1}{w} - wq \big) & 0 & \frac{1}{w} + wq \vspace{1mm}\\
 0 & - \big(\frac{1}{w} + wq \big) & 0 \end{pmatrix} \qquad \text{and} \qquad
 \widehat{N}_3 = wr \begin{pmatrix}
 \hphantom{-}0 & {\rm i} & \hphantom{-}0 \\ -{\rm i} & 0 & -1 \\ \hphantom{-}0 & 1 & \hphantom{-}0
 \end{pmatrix}\!,
\end{gather*}
is an $\mathfrak{so}(3, C_K)$ system corresponding to the linear differential equation~\eqref{mdar1}. A fundamental matrix for this system is
 \begin{gather}\label{eq:sol-fund-so3-2}
 \mathbf{Z}_1 = S \cdot \mathbf{Y}_1 = \begin{pmatrix}
 y_1^2 +w^2 (y_1')^2 & y_1 y_2 + w^2 y_1' y_2' & y_2^2 + w^2 (y_2')^2 \vspace{1mm}\\
 2{\rm i}w y_1 y_1' & {\rm i}w (y_1 y_2' + y_1' y_2) & 2{\rm i}w y_2 y_2' \vspace{1mm}\\
 {\rm i}\big(y_1^2 - w^2 (y_1')^2\big) & {\rm i}\big(y_1 y_2 - w^2 y_1' y_2'\big) & {\rm i}\big(y_2^2 - w^2 (y_2')^2\big)
\end{pmatrix}\!,
 \end{gather}
 where $\mathbf{Y}_1$ is given by~\eqref{eq:sol-fund-sl-2} and $\{ y_1, y_2 \}$ is a basis of solutions of equation~\eqref{eqdarb}.

 Then, the second expression for the Darboux transformation for $\mathfrak{so}(3, C_K)$ systems is given by the next corollary.

\begin{Corollary}[second Darboux transformation for $\mathfrak{so}(3, C_K)$]\label{secthdars032}
Consider the system $\big[\widehat{\Omega}_m\big]$\emph{:} $Z_1' = - \big(\widehat{\Omega}_0 + m \widehat{N}_3 \big) Z_1$ given by~\eqref{eq:so-system-2}. Let $T_{2,m}$ be the matrix
\begin{align*}
 T_{2,m} &= S P_{2,m} S^{-1} \nonumber
 \\
 &= \dfrac{1}{2r} \begin{pmatrix} {w}^{2}+{\rho}^{2}+{\theta_0 }^{2}+{\frac {{\nu}^{2}}{{w}^{2}}} & 2{\rm i} \big (w \theta_0 -{\frac {\nu \rho}{w}} \big )
& {\rm i} \big ( {w}^{2}+ {\rho}^{2}- {\theta_0 }^{2}-{\frac {{\nu}^{2}}{{w}^{2}}} \big )
\\[1mm]
2 {\rm i} \big ( w \rho-{\frac {\nu \theta_0 }{w}} \big ) & 2 ( \nu - \rho \theta_0 ) & -2 \big ( w \rho + {\frac { \nu \theta_0 }{w}} \big )
\\[1mm]
{\rm i} \big ( {w}^{2}- {\rho}^{2}+ {\theta_0 }^{2}-{\frac {{\nu}^{2}}{{w}^{2}}} \big ) & -2 \big( w \theta_0 + {\frac {\nu \rho}{w}} \big ) & -{
w}^{2}+{\rho}^{2}+{\theta_0 }^{2}-{\frac {{\nu}^{2}}{{w}^{2}}}
\end{pmatrix}\!,
\end{align*}
where matrix $P_{2,m}$ is defined by expression~\eqref{eq:DT sl22} and $\nu=mr-\theta_0 \rho$. Then, the gauge transformation $T_{2,m}$ is a Darboux transformation for system $\big[\widehat{\Omega}_m\big]$.
\end{Corollary}

\begin{proof}
Since $ \widetilde{Y}_1 =P_{2,m} \cdot Y_1$ by Corollary~\ref{thnudarsl22},
it follows that
\begin{gather*}
T_{2,m} = S \cdot P_{2,m} \cdot S^{-1}.
\end{gather*}
By construction, the image of $\widehat{\Omega}_0 + m \widehat{N}_3$ by this gauge transformation is
$\widetilde{\widehat{\Omega}}_0 + m \widehat{N}_3$, where~$\widetilde{\widehat{\Omega}}_0$ is obtained from $\widehat{\Omega}_0$ by changing $q$ by the function $q_{1}$ obtained in~\eqref{etda13} by the Darboux transformation:
\begin{gather*}
\widetilde{\widehat{\Omega}}_0 = \begin{pmatrix}
 0 & {\rm i} \big (\frac{1}{w} - w\tilde{q} \big) & 0
 \vspace{1mm}\\
 -{\rm i} \big (\frac{1}{w} - w\tilde{q} \big) & 0 & \frac{1}{w} + w\tilde{q} \\[5pt] 0 & - \big ( \frac{1}{w} + w\tilde{q} \big ) & 0 \end{pmatrix}\!.
 \end{gather*}
The rest of the corollary is proved following the same argument as in Proposition~\ref{prop1}.
\end{proof}

As in all the above results, the transformation can be factored into an $m$-dependent part and an independent part. The formula is not as compact as the previous ones but is easily found with a computer algebra system once one is equipped with this paper's methodology and results.

We note that Darboux transformations $T_{1,m}$ and $T_{2,m}$ are not equivalent when $w\neq1$ because matrices $Q$ and $S$ are then different.

Corollaries~\ref{thnudarsl22} and~\ref{secthdars032} can be summarized in the diagram below. Each arrow in this diagram is now made explicit and this provides a complete algorithmic way to produce Darboux transformations for any family of orthogonal systems:
\begin{gather}\label{eq:diag-2}
\begin{split}
 \xymatrix{ [A_m]\colon X'= -(A_0 +mN) X \ar[rr]^{P_m} \ar@<1ex>[d]^{\Delta} & & \big[\widetilde{A}_m\big]\colon \widetilde{X}'= -\big(\widetilde{A}_0 +mN\big) \widetilde{X} \ar@<1ex>[d]^{\Delta}
 \\
 [B_m]\colon X_1'= -(B_0 + mN_1) X_1 \ar[rr]^{\Delta P_m \Delta^{-1}} \ar@<1ex>[d]^{\sym^2} \ar[u] & & \big[\widetilde{B}_{m}\big]\colon \widetilde{X}_{1}'=-\big(\widetilde{B}_0+mN_1\big) \widetilde{X}_{1} \ar@<1ex>[d]^{\sym^2} \ar[u]
 \\
 \big[\sym^2 (B_m)\big]\colon Y_1'=-\big(\widehat{S}_2 + m \widehat{N}_2\big) Y_1 \ar[rr]^{P_{2,m}} \ar@<1ex>[d]^{S} \ar[u] & & \big[\sym^2 \big(\widetilde{B}_m\big)\big]\colon \widetilde{Y}_1'=-\big(\widetilde{\widehat{S}}_2 + m \widehat{N}_2\big) \widetilde{Y}_1 \ar@<1ex>[d]^{S} \ar[u]
 \\
 \big[\widehat{\Omega}_m\big]\colon Z_1' =-\big(\widehat{\Omega}_0 + m \widehat{N}_3\big) Z_1 \ar[rr]^{T_{2,m}} \ar[u] & & \big[\widetilde{\widehat{\Omega}}_m\big]\colon \widetilde{Z}_1' =-\big(\widetilde{\widehat{\Omega}}_0 + m \widehat{N}_3\big) \widetilde{Z}_1. \ar[u] }
 \end{split}\hspace{-10mm}
\end{gather}

\subsection{Extension to general differential systems with an orthogonal Galois group}
\label{extension}

We now show how all the above results can be extended to the general case of a linear differential system $Y' = -A Y$ whose differential Galois group is in the special orthogonal group ${\rm SO}(3, C_K)$, even when the
 matrix $A$ is not in $\mathfrak{so}(3, C_K)$. This is important as many linear differential systems arising from physical models come equipped with this orthogonality property without having their matrix in $\mathfrak{so}(3, C_K)$. The presence of an orthogonal Galois group means that there is a gauge equivalence between the system and its adjoint, a duality; in physics, this is often inherited from a duality property. This has been proved by Bogner for Calabi--Yau operators, see~\cite{Bo13a,BoRe13a}, and studied in a slightly more general context in~\cite{BoBoMaWe15a,BoHaMaWe14a,BoHaMaWe15a} for operators arising from statistical mechanics or combinatorics.

First we recall how to detect this situation, namely how to decide
when the differential Galois group is in ${\rm SO}(3, C_K)$, see, e.g.,~\cite{NgPu10a,Si88a,Ho07a} or the old book of Darboux~\cite[pp.~28--29]{da3}. The system should be irreducible (it has no hyperexponential solution), the trace of $A$ should be a~logarithmic derivative (i.e., the equation $y'=-\operatorname{Tr}(A)y$ has a rational solution) and the second symmetric power system $Z' = -\sym^2(A) Z$ should have a rational solution, corresponding to the quadratic invariant of the special orthogonal group. All these properties can be effectively tested in a computer algebra system, see~\cite{BaClElWe12a} and references therein (notably an implementation in the computer algebra system \textsc{Maple}).
When they are fulfilled, the differential Galois group is in ${\rm SO}(3, C_K)$.

In general, the matrix of such a differential systems is not in $\mathfrak{so}(3, C_K)$; then, the results of the previous section would not apply directly. However, using the constructive theory of reduced forms from~\cite{BaClDiWe20a, BaClDiWe16a}, we may compute effectively a gauge transformation matrix $P$ such that, letting $Y=PZ$, the new unknown $Z$ satisfies a system $Z'=-BZ$ with $B\in \mathfrak{so}(3, K)$. Then the results of the previous sections apply: we may solve using solutions of second-order equations and construct families of equations of similar shapes via Darboux transformation.
Earlier versions of such a ``reduction'' process also appears in Singer's work~\cite{Si88a} and several subsequent works on solving linear differential equations in terms of lower-order equations, notably~\cite{Ho07a} (which comes with an efficient \textsc{Maple} implementation) and~\cite{NgPu10a} which is more general.

Note that these observations already appears in the old book of Darboux~\cite[p.~28--29]{da3}: he shows (in old language) how to identify third-order linear differential systems with an orthogonal Galois group using a first integral (a quadratic invariant, in our language); he then shows how to transform such a system into the orthogonal form treated in this paper and calls it the \emph{type} or 	\emph{``la forme r\'eduite''} of the class of third-order systems admitting a quadratic first integral.

In conclusion: given a third-order linear differential system, we can check algorithmically whether it has a Galois group in ${\rm SO}(3)$. When this is the case, we can find a gauge transformation which reduces the systems into an orthogonal one. Then, we may apply the machinery of the previous parts to obtain Darboux transformations. This shows that the machinery of this paper allows the reader to construct Darboux transformations for any third-order linear differential system with an orthogonal Galois group or, equivalently, with a quadratic first integral.

\section{Applications}\label{applications}
In this section, to motivate the results of this paper, we present some examples coming from supersymmetric quantum mechanics and differential geometry.

\subsection{Supersymmetric quantum mechanics}\label{sec:SUSY}

The Schr\"odinger equation for the stationary and non-relativistic case is given by
\begin{gather*}
H\psi=\lambda \psi,\qquad H=-\partial_x^2 +V(x),
\end{gather*}
where $\lambda$ is called the \emph{energy}, $V$ is called the \emph{potential} and $\psi$ is called the \emph{wave function}. Supersymmetric quantum mechanics in the Witten's formalism was introduced by himself in~\cite[Section~6]{wi} as a toy model. Witten introduced the \emph{supersymmetric partner Hamiltonians} $H_\pm$ as follows
\begin{gather*}
H_\pm=- \partial_x^2 +V_\pm(x),\qquad V_\pm=W^2\pm W',
\end{gather*}
where $V_\pm$ are called the \emph{supersymmetric partner potentials} and $W$ is called the \emph{superpotential} which satisfies
\begin{gather*}
W=-\dfrac{\psi_0'}{\psi_0},\qquad H_-\psi_0=\lambda_0\psi_0,
\end{gather*}
where $\psi_0$ is called the \emph{ground state} and $\lambda_0$ is an specific value of the energy $\lambda$.

We can go from $H_-\psi=\lambda\psi$ to $H_+\widetilde{\psi}=\lambda\widetilde{\psi}$ through a Darboux transformation, where $\theta=\psi_0$, $m=-\lambda$, $y=\psi$, $u=\widetilde{\psi}$, $p=0$, $q=-V$, and $r=1$.

Gendenshtein in~\cite{ge} introduced what today is called \emph{shape invariant} potentials, that is, potentials with the \emph{shape invariance property}: the potential $V=V_-=V_-(x;a)$ has the shape invariance property if and only if its supersymmetric partner potential $V_+=V_+(x;a)$ can be written as $V_+(x;a)=V_-(x;a_1)+R(a_1)$, where $a$ is a set of
parameters, $a_1$ is a new set of parameters given by $a_1=f(a)$ for some function $f$ and $R$ is a remainder which does not depend on $x$, see~\cite[Introduction and Section~4.3]{AP1}.
In other words, the supersymmetric partner potentials differs only in parameters. Applying systematically this procedure, one gets
iterated values $a_k$ of parameters and can obtain
the spectrum as values of energy $\lambda$ given by
\begin{gather*}
\lambda=\sum_{k=1}^n R(a_k).
\end{gather*}
See~\cite[Section~4.3]{AP1} for explanations on this phenomenon and how to detect it.

 Moreover, $H_-=A^\dag A$ and $H_+=A A^\dag$, where $A^\dag=-\partial_x+W$ and $A=\partial_x+W$ are called the \emph{ladder $($raising and lowering$)$ operators}, see~\cite{DKS}. We can rewrite the starting potential $V$ as ${V_--\lambda_0}$ to apply Darboux transformations. Thus, we can obtain the rest of wave functions applying it as follows: $\psi_1=A(x;a_0)^\dag\psi_0$, and in general as $\psi_k=A(x;a_{k-1})^\dag\psi_{k-1}$, where ${a_0=a}$, $a_1=f(a_0)$.
First examples of rational shape invariant potentials correspond to harmonic oscillator and Coulomb potentials, for one dimensional and three-dimensional cases.

In the following, we combine this theoretical background regarding the Schr\"odinger equation and supersymmetric quantum mechanics with our results from the previous section.

\subsubsection[Darboux transformation in matrix form for supersymmetric quantum mechanics]
{Darboux transformation in matrix form for supersymmetric \\quantum mechanics}
We apply our previous results to a {matrix form} of the Schr\"odinger equation. We start by introducing the following $2\times 2$ matrix Schr\"odinger operators, with supersymmetric partner potentials, related to the systems~\eqref{mdar1} and~\eqref{mdar2} as follows:
\begin{gather*}
\mathcal H_\pm=-\partial_x+\mathbf{V}_\pm,\qquad
\mathbf{V}_{\pm} =\begin{pmatrix}0&1\\V_{\pm}&0\end{pmatrix}\!,\qquad
\mathbf{V}_{-}=-A_0,\qquad
\mathbf{V}_{+}=-\widetilde{A}_0,
\end{gather*}
where
\begin{gather*}
\mathcal H_-\Psi=\mathfrak{E}_{\lambda} \Psi,\qquad
H_+\widetilde{\Psi}=\mathfrak{E}_{\lambda} \widetilde{\Psi}, \qquad
\Psi=\begin{pmatrix} \psi\\ \psi'\end{pmatrix}\!,\qquad
\widetilde{\Psi}=\begin{pmatrix} \widetilde{\psi}\\ \widetilde{\psi'}\end{pmatrix}\!,\\ \mathfrak{E}_{\lambda}=-\lambda N,\qquad
-N=\begin{pmatrix}0&0\\1&0\end{pmatrix}\!.
\end{gather*}
According to Proposition~\ref{prop1}, the relevant Darboux transformation in this $2\times 2$ matrix formalism is given by
\begin{gather*}
\widetilde{\Psi}=P_\lambda\Psi,	\qquad
P_\lambda=\begin{pmatrix}W&1\\W^2-\lambda&W \end{pmatrix}=
\begin{pmatrix} 0 & 1 \\ -\lambda&W\end{pmatrix}\cdot \begin{pmatrix} 1&0\\W&1\end{pmatrix}\!,
\end{gather*}
and the supersymmetric partner potentials $\mathbf{V}_\pm$ will depend on the supersymmetric partner potentials $V_\pm$ according to the original Witten's formalism, i.e., $V_+=V_-+2W'$, which lead us to $\mathbf{V}_+=\mathbf{V}_--2W'N$.

Now we present the shape invariance property for this $2\times 2$ matrix formalism of Schr\"odinger equation as follows. Consider the parametric supersymmetric partner potentials and matrix $\mathbf{R}(a_1)$ as follows:
\begin{gather*}
\mathbf{V}_\pm(x;a)=\begin{pmatrix}0&1\\V_{\pm}(x;a)&0\end{pmatrix}\!,\qquad
\mathbf{R}(a_1)=-R(a_1)N,
\end{gather*}where, as in the classical case, $a$ is a set of parameters and $a_1=f(a)$. The potential $\mathbf{V}=\mathbf{V}_-=\mathbf{V}_-(x;a)$ has the shape invariance property if and only if its supersymmetric partner potential satisfies
\begin{gather*}
\mathbf{V}_+=\mathbf{V}_+(x;a)=\mathbf{V}_-(x;a_1)+\mathbf{R}(a_1).
\end{gather*}
We can see again that supersymmetric partner potentials differ only in parameters. Applying systematically this procedure, we can obtain the spectrum as the values of energy $\mathfrak E_\lambda (a)$, where $ \mathfrak E_\lambda (1) = \mathfrak E_\lambda $, such that
\begin{gather*}
\mathfrak E_\lambda (a)=\sum_{k=1}^n \mathbf{R}(a_k).
\end{gather*}
Also we present the ladder operators for this $2\times 2$ matrix formalism of Schr\"odinger equation:
\begin{gather*}
\mathbf{A}=\begin{pmatrix}A&0\\WA&0\end{pmatrix}\!,\qquad
\mathbf{A}^\dag=\begin{pmatrix} A^\dag &0\\2W'-WA^\dag &0 \end{pmatrix}\!.
\end{gather*}

We illustrate this formalism with the 1D-harmonic oscillator, which is a classical rational shape invariant potential. The superpotential for harmonic oscillator is $W=x$, thus the supersymmetric partner potentials are given by
\begin{gather*}
\mathbf{V}_-=\begin{pmatrix}0&1\\x^2- 1&0\end{pmatrix}\!, \qquad \mathbf{V}_+=\begin{pmatrix}0&1\\x^2+1&0\end{pmatrix}=\mathbf{V}_-+\begin{pmatrix}0&0\\2&0\end{pmatrix}\!.
\end{gather*}
Therefore, introducing a multiplicative parameter $a$ in $\mathbf{V}_-$, such that $\mathbf{V}_-(x;a)=\mathbf{V}_-(x)$ for $a=1$, we obtain
\begin{gather*}
f(a_1)=2a,\qquad \mathbf{R}(a_1)=\begin{pmatrix}0&0\\2a&0\end{pmatrix}\!,\qquad \mathfrak{E}_{\lambda} (a) = \sum_{k=1}^n \mathbf{R}(a_k)=\begin{pmatrix}0&0\\2na&0\end{pmatrix}\!.
\end{gather*}
Thus, for $a=1$ the spectrum of $\mathcal H_-$ is
\begin{gather*}
\mathrm{Spec}(\mathcal H_-)=\{\mathfrak E_\lambda\colon\lambda\in 2\mathbb{Z}_+\}.
\end{gather*}
For instance, we have
\begin{gather*}
\widetilde{\Psi}=P_\lambda\Psi=\begin{pmatrix}
x&1\\x^2-\lambda&x\end{pmatrix}\begin{pmatrix}
H_{\frac{\lambda}{2}}(x)\\H'_{\frac{\lambda}{2}}(x)-xH_{\frac{\lambda}{2}}(x)
\end{pmatrix}\exp\left({-}\frac{x^2}{2}\right),
\end{gather*}
where $H_{\frac{\lambda}{2}}$ denotes the Hermite polynomial of degree $\frac{\lambda}{2}$. The ladder operators are given respectively by
\begin{gather*}
\mathbf{A}=\begin{pmatrix}
A&0\\xA&0\end{pmatrix}\!,\qquad \mathbf{A}^\dag=\begin{pmatrix}
A^\dag&0\\2-xA^\dag&0\end{pmatrix}\!,
\end{gather*}
where $A= \partial_x + x$ and $A^\dag = -\partial_x + x$. Using these ladder operators we can obtain
\begin{gather*}
\Psi_n=\mathbf{A}^\dagger \Psi_{n-1}(x), \qquad \textrm{where} \quad
\lambda=2n \quad \textrm{and} \quad \Psi_0=\begin{pmatrix}\exp\bigl(-\frac{x^2}{2}\bigr)\\[1mm] -x\exp\bigl(-\frac{x^2}{2}\bigr)\end{pmatrix}\!.
\end{gather*}

\subsubsection{Second symmetric power approach for supersymmetric quantum mechanics}

The following results correspond to the second symmetric power of Schr\"odinger equation in the previous matrix formalism. Thus, we obtain the following $3\times 3$ matrix Schr\"odinger operators, with supersymmetric partner potentials according to equations~\eqref{ssyp m},~\eqref{ssyp m dt},~\eqref{eq:sist ssmsl2} and~\eqref{ssyp sl2} as follows:
\begin{gather*}
\mathcal H_\pm=-\partial_x+\mathbf{V}_\pm,\qquad \mathbf{V}_-=-S=-\widehat{S},\qquad \mathbf{V}_+=-\widetilde{S}=-\widetilde{\widehat{S}},\qquad \mathbf{V}_{\pm} =\begin{pmatrix}0&1&0\\2V_{\pm}&0&2\\0&V_\pm&0\end{pmatrix}\!,
\end{gather*}
where
\begin{gather*}
\mathcal H_-\Psi=\mathfrak{E}_{\lambda} \Psi,\qquad
H_+\widetilde{\Psi}=\mathfrak{E}_{\lambda} \widetilde{\Psi}, \qquad
\Psi=\begin{pmatrix} \psi^2\\2\psi\psi'\\(\psi')^2\end{pmatrix}\!,\qquad \widetilde{\Psi}=\begin{pmatrix} \widetilde{\psi}^2\\2\widetilde{\psi}\widetilde{\psi'}\\(\widetilde{\psi}')^2\end{pmatrix}\!,\qquad \mathfrak{E}_{\lambda}=-\lambda N_1,
\end{gather*}
and
\begin{gather*}
-N_1=-N_2=\begin{pmatrix}0&0&0\\2&0&0\\0&1&0\end{pmatrix}\!.
\end{gather*}
Using~\eqref{eq:DT sl2}, our generalized Darboux transformation in this $3\times 3$ matrix formalism is given by
\begin{gather*}
\widetilde{\Psi}=P_\lambda\Psi,\qquad P_\lambda=P_{1,\lambda}=P_{2,\lambda}=\begin{pmatrix}
W^2&W&1\\ 2W\big(W^2-\lambda\big)&2W^2-\lambda&2W\\ \big(W^2-\lambda\big)^2&W\big(W^2-\lambda\big)&W^2
\end{pmatrix}
\end{gather*}
with the factorization into a $\lambda$-dependent part and a part with only $W$:
\begin{gather*}
P_\lambda =
\begin{pmatrix} 0 & 0 & 1 \\ 0 & - \lambda & 2 W \\ \lambda^2 & -\lambda W & W^2 \end{pmatrix} \cdot \begin{pmatrix} 1 & 0 & 0 \\ 2W & 1 & 0 \\ W^2 & W & 1 \end{pmatrix}\!.
\end{gather*}
As in the previous case, the supersymmetric partner potentials $\mathbf{V}_\pm$ will depend on the supersymmetric partner potentials $V_\pm$ according to the original Witten's formalism, i.e., $V_+=V_-+2W'$, which lead us to
\begin{gather*}
\mathbf{V}_+=\mathbf{V}_- -2W'N_1,
\end{gather*}
Now we present the shape invariance property for this $3\times 3$ matrix formalism of Schr\"odinger equation as follows. Consider the parametric supersymmetric partner potentials and matrix~$\mathbf{R}(a_1)$ as follows:
\begin{gather*}
\mathbf{V}_\pm(x;a)=\begin{pmatrix}0&1&0\\2V_{\pm}(x;a)&0&2\\0&V_{\pm}(x;a)&0\end{pmatrix}\!,\qquad \mathbf{R}(a_1)=-R(a_1)N,
\end{gather*}
where, as in the previous case, $a$ is a set of parameters and $a_1=f(a)$. The potential $\mathbf{V}=\mathbf{V}_-=\mathbf{V}_-(x;a)$ has the shape invariance property if and only if its supersymmetric partner potential can be written as
\begin{gather*}
\mathbf{V}_+=\mathbf{V}_+(x;a)=\mathbf{V}_-(x;a_1)+\mathbf{R}(a_1).
\end{gather*}
We can see again that supersymmetric partner potentials differ only in parameters. Applying systematically this procedure, we can obtain the spectrum as the values of energy $\mathfrak E_\lambda (a)$
given by $ \mathfrak E_\lambda (1) = \mathfrak E_\lambda $ and
\begin{gather*}
\mathfrak E_\lambda (a)=\sum_{k=1}^n \mathbf{R}(a_k).
\end{gather*}
The ladder operators for this $3\times 3$ matrix formalism of Schr\"odinger equation are
\begin{gather*}
\mathbf{A}=\begin{pmatrix}WA&0&1\\2W^2A&0&2W\\W^3A&0&W^2\end{pmatrix}\!, \qquad
\mathbf{A}^\dag=\begin{pmatrix}
WA^\dag&0&1\\-2W^2A^\dag&0&-2W\\W^3A^\dag&0&W^2
\end{pmatrix}\!.
\end{gather*}

We illustrate this formalism with the 1D-harmonic oscillator, which is a classical rational shape invariant potential. The superpotential for harmonic oscillator is $W=x$, thus the supersymmetric partner potentials are given by
\begin{gather*}
\mathbf{V}_-=\begin{pmatrix}0&1&0\\2x^2- 2&0&2\\0&x^2-1&0\end{pmatrix}\!, \qquad \mathbf{V}_+=\begin{pmatrix}0&1&0\\2x^2+2&0&2\\0&x^2+1&0\end{pmatrix}
=\mathbf{V}_-+\begin{pmatrix}0&0&0\\4&0&0\\0&2&0\end{pmatrix}\!.
\end{gather*}
Therefore, introducing a multiplicative parameter $a$ in $\mathbf{V}_-$, such that $\mathbf{V}_-(x;a)=\mathbf{V}_-(x)$ for $a=1$, we obtain
\begin{gather*}
f(a_1)=2a,\qquad \mathbf{R}(a_1)=\begin{pmatrix}0&0&0\\4a&0&0\\0&2a&0\end{pmatrix}\!,\qquad \mathfrak{E}_{\lambda} (a) = \sum_{k=1}^n \mathbf{R}(a_k)=\begin{pmatrix}0&0&0\\4na&0&0\\0&2na&0\end{pmatrix}\!.
\end{gather*}
Thus, for $a=1$ the spectrum of $\mathcal H_-$ is
\begin{gather*}
\mathrm{Spec}(\mathcal H_-)=\{\mathfrak E_\lambda\colon\lambda\in 2\mathbb{Z}_+\}.
\end{gather*}
For instance, we have
\begin{gather*}
\widetilde{\Psi}=P_\lambda\Psi=\begin{pmatrix}
x^2&x&1\\2x^3-2\lambda x&2x^2-\lambda&2x\\x^4-2\lambda x^2+\lambda^2&x^3-\lambda x&x^2
\end{pmatrix}\begin{pmatrix}
H^2_{\frac{\lambda}{2}}(x)
\\
\big(H^2\big)'_{\frac{\lambda}{2}}(x)-2xH^2_{\frac{\lambda}{2}}(x)
\\
\big(H'_{\frac{\lambda}{2}}(x)-xH_{\frac{\lambda}{2}}(x)\big)^2
\end{pmatrix}\exp\bigl(-x^2\bigr),
\end{gather*}
where $H_{\frac{\lambda}{2}}$ denotes the Hermite polynomial of degree $\frac{\lambda}{2}$. The ladder operators are given respectively by
\begin{gather*}
\mathbf{A}=\begin{pmatrix}
xA&0&1\\2x^2A&0&2x\\x^3A&0&x^2
\end{pmatrix}\!, \qquad
\mathbf{A}^\dag=\begin{pmatrix}
xA^\dag&0&1\\-2x^2A^\dag&0&-2x\\x^3A^\dag&0&x^2
\end{pmatrix}\!.
\end{gather*}

\subsection{Some \texorpdfstring{\mathversion{bold}$\mathfrak{so}(3, C_K)$}{so(3,C\_K)} systems}\label{sec:applic-so3}

In this section, we revisit from the point of view developed in this article two well-known problems which arise expressed as $\mathfrak{so}(3, C_K)$ systems.

\subsubsection{Frenet--Serret formulas}

Given a nondegenerate curve in the space, denote by $T$ the tangent unit vector to the curve, by~$N$ the normal unit vector and by $B = T \times N$ the binormal unit vector. Then, the Frenet--Serret formulas can be formulated as the following $\mathfrak{so}(3, C_K)$ system:
\begin{gather}
\begin{pmatrix} T \\ N \\ B \end{pmatrix}' = - \begin{pmatrix}
0 & -\kappa & 0 \\ \kappa & 0 & - \tau \\ 0 & \tau & 0 \end{pmatrix} \cdot
\begin{pmatrix} T \\ N \\ B \end{pmatrix} = - \Omega \cdot \begin{pmatrix} T \\ N \\ B \end{pmatrix}\!,
\end{gather}
where $'$ denotes the derivative with respect to arclength, $\kappa$ is the curvature of the curve and $\tau$ is its torsion, see~\cite[Chapter~1]{MSp} for more details.

In order to apply our previous formalism to this $\mathfrak{so}(3, C_K)$ system, we have two possibilities: we can use system~\eqref{eq:system so3} or system~\eqref{eq:so-system-2}.

In the first case, the identification of matrix $\Omega$ with matrix $\Omega_0$ given by Lemma~\ref{gaugemat} yields the degenerate situation: $\kappa = -{\rm i}p $, $\tau = {\rm i} (q-1)$ and $0 = q+1$, hence, $q= -1$ and $\tau= -2{\rm i}$. The second-order linear differential equation associated to this system is
\begin{gather}\label{eq:FS-eq1}
 \mathcal{L} y = y'' +{\rm i}\kappa y' -y =0.
\end{gather}
Since $p= \frac{w'}{w} ={\rm i} \kappa$, we find that $w = {\rm e}^{({\rm i} \int \kappa {\rm d}x)}$. As immediate application of $\mathfrak{so}(3, C_K)$ systems to Frenet--Serret formulas, we obtain directly a fundamental matrix of solutions through Corollary~\ref{corso3}:
\begin{gather*}
 \mathbf{Z} = {\rm e}^{({\rm i} \int \kappa dx)} \cdot
\begin{pmatrix}
y_1^2 - (y_1')^2 & y_1 y_2 - y_1' y_2' & y_2^2 -(y_2')^{2} \\[3pt]
{\rm i} \big(y_1^2 + (y_1')^{2}\big) & {\rm i} (y_1 y_2 + y_1' y_2') & {\rm i} \big(y_2^2 + (y_2')^{2}\big) \\[3pt]
-2 y_1 y_1' & - y_1 y_2' - y_1' y_2 & -2 y_2 y_2' \end{pmatrix}\!,
\end{gather*}
where $\{ y_1, y_2\}$ is a basis of solutions of equation~\eqref{eq:FS-eq1}.

Notice that this framework is only valid for curves with torsion $\tau=-2{\rm i}$. However, we can avoid this restriction by using the second approach, given by equation~\eqref{eq:so-system-2}. In this case, the gauge change given by $S$ transform the matrix $\Omega_0$ into a matrix with the same structure as $\Omega$, namely $\widehat{\Omega}_0$. Identifying the entries of both matrices we obtain: $\kappa = -{\rm i}(1/w -wq)$ and $\tau = -(1/w +wq)$. Thus, $w = \frac{2}{{\rm i}\kappa - \tau}$ and $q = \frac{\kappa^2 + \tau^2}{4}$ and the second-order linear differential equation associated to the system in this case is
\begin{gather}\label{eq:FS-eq2}
 \mathcal{L} y = y'' - \dfrac{({\rm i}\kappa - \tau)'}{{\rm i} \kappa - \tau} y' +\dfrac{\kappa^2 + \tau^2}{4}y =0.
\end{gather}
Under this assumptions, a fundamental matrix for this system, given by~\eqref{eq:sol-fund-so3-2}, becomes{\samepage
\begin{gather*}
\mathbf{Z}_1 = \begin{pmatrix}
y_1^2 + \dfrac{4(y_1')^2}{({\rm i}\kappa - \tau)^2} & y_1 y_2 + \dfrac{4y_1' y_2'}{({\rm i}\kappa - \tau)^2} & y_2^2 + \dfrac{4 (y_2')^2}{({\rm i}\kappa - \tau)^2} \\[9pt]
\dfrac{4{\rm i} y_1 y_1'}{{\rm i}\kappa - \tau} & \dfrac{2{\rm i} (y_1 y_2' + y_1' y_2)}{{\rm i}\kappa - \tau} & \dfrac{4{\rm i} y_2 y_2'}{{\rm i}\kappa - \tau} \\[9pt]
{\rm i} \bigg( y_1^2 - \dfrac{4 (y_1')^2}{({\rm i}\kappa - \tau)^2}\bigg) & {\rm i} \bigg( y_1 y_2 - \dfrac{4 y_1' y_2'}{({\rm i}\kappa - \tau)^2} \bigg) & {\rm i} \bigg( y_2^2 - \dfrac{4 (y_2')^2}{({\rm i}\kappa - \tau)^2} \bigg ) \end{pmatrix}\!,
\end{gather*}
where $\{ y_1, y_2\}$ form a basis of solutions of the second-order equation~\eqref{eq:FS-eq2}.}

Finally, we consider the following
variant of the Frenet--Serret system using the second approach for $\mathfrak{so}(3, C_K)$ systems:
\begin{gather*}
\begin{pmatrix} T \\ N \\ B \end{pmatrix}' = - \left ( \begin{pmatrix}
0 & -\kappa & 0 \\ \kappa & 0 & - \tau \\ 0 & \tau & 0 \end{pmatrix} + \frac{2m}{{\rm i}\kappa - \tau} \begin{pmatrix} 0 & {\rm i} & 0 \\ -{\rm i} & 0 & -1 \\ 0 & 1 & 0 \end{pmatrix} \right ) \cdot \begin{pmatrix} T \\ N \\ B \end{pmatrix} .
\end{gather*}
This equation is obtained by directly applying equation~\eqref{eq:so-system-2}
 for $r=1$ and the expressions for~$w$ and $q$ obtained before.

Following the philosophy of Darboux we can construct an infinite chain of such perturbed Frenet--Serret systems applying Corollary~\ref{secthdars032} (we restrict ourselves to the case $r=1$ since it is the usual situation in the applications, see, for instance,
the classical book of Matveev and Salle~\cite{MaSa91a}). The Darboux transformation is given by
\begin{gather*}
T_{2,m} = \dfrac{1}{2} \begin{pmatrix} \dfrac{4}{\eta^2} + \rho^2 + \theta_0 ^2 + \dfrac{\nu^2 \eta^2}{4} & {\rm i} \bigg (\dfrac{4\theta_0 }{\eta} - \nu \rho \eta \bigg ) & {\rm i} \bigg ( \dfrac{4}{\eta^2} + \rho^2 - \theta_0 ^2 - \dfrac{\nu^2 \eta^2}{4} \bigg )\\[10pt]
{\rm i} \bigg ( \dfrac{4\rho}{\eta} - \nu \theta_0 \eta \bigg ) & 2 ( \nu -\rho \theta_0 ) & \dfrac{-4\rho}{\eta} - \nu \theta_0 \eta \\[10pt]
{\rm i} \bigg ( \dfrac{4}{\eta^2} - \rho^2 + \theta_0 ^2 - \dfrac{\nu^2 \eta^2}{4} \bigg ) & \dfrac{-4\theta_0 }{\eta} - \nu \rho \eta & \rho^2 + \theta_0 ^2 - \dfrac{4}{\eta^2} - \dfrac{\nu^2 \eta^2}{4}
\end{pmatrix}
\end{gather*}
for $\theta_0 = \frac{y'}{y}$, $\rho = -\theta_0 + \frac{\eta'}{\eta}$, $\nu=m -\theta_0 \rho$ and $\eta = i\kappa - \tau$, where $y$ is a solution of equation~\eqref{eq:FS-eq2}.
As before, this matrix could be factored into an $m$-dependent and an independent part.

\subsubsection{Rigid solid problem}

A rigid solid consists of a set of points in the space preserving the distance among them under the action of some applied forces. The transformations allowed for a rigid solid are translations and rotations.

The Poisson equation describes the motion of the rigid body in space:
\begin{gather*}
\gamma' = \gamma \times \omega ,
\end{gather*}
where $\gamma = (\gamma_1, \gamma_2, \gamma_3)^{\rm T}$ is a unit vector fixed in space, $\omega = (\omega_1, \omega_2, \omega_3)^{\rm T}$ the angular velocity vector and ${}' = \partial_t$. See~\cite{acjios, femapr, femapr0} for more details.

We follow the constraints and notation of~\cite{femapr}. Hence, we take $\omega_3 =0$ and restrict ourselves to rigid transformations in the plane. In this case, the Poisson equation can be rewritten as the $\mathfrak{so}(3, C_K)$ system:
\begin{gather*}
 \begin{pmatrix} \gamma_1 \\ \gamma_2 \\ \gamma_3 \end{pmatrix}' = - \begin{pmatrix}
0 & 0 & \omega_2 \\ 0 & 0 & - \omega_1 \\ -\omega_2 & \omega_1 & 0 \end{pmatrix} \cdot \begin{pmatrix} \gamma_1 \\ \gamma_2 \\ \gamma_3 \end{pmatrix} = - A \cdot \begin{pmatrix} \gamma_1 \\ \gamma_2 \\ \gamma_3 \end{pmatrix}\!.
\end{gather*}

Fedorov et al.\ studied in~\cite{femapr} this system for a general case of matrix~$A$. Acosta-Hum\'anez et al.\ considered in~\cite{acjios} a particular case with $\omega_1 = \frac{{\rm e}^x - {\rm e}^{-x}}{{\rm e}^x + {\rm e}^{-x}}$ and $\omega_2 = \frac{2 \sqrt{2}}{{\rm e}^x + {\rm e}^{-x}}$. In this work, we are going to restrict the rigid transformations allowed to the coupled case ${\rm i}\omega_1 + \omega_2 = 2$.
For that, we consider the $\mathfrak{so}(3, C_K)$ system given in Lemma~\ref{gaugemat} and apply the formalism developed for the $\mathfrak{so}(3, C_K)$ system~\eqref{eq:system so3}.
The identification of matrix $A$ with matrix $\Omega$ leads to: $p=0$, $\omega_1 = {\rm i}(q-1)$, $\omega_2 = q+1$, hence, $q = 1- {\rm i}\omega_1 = \omega_2 -1$. This yields the announced coupled situation $i\omega_1 + \omega_2 = 2$.

The second-order linear differential equation associated to this problem is
\begin{gather}\label{eq:RS-eq1}
 \mathcal{L} y = y'' + (1- i\omega_1)y = y'' + (\omega_2 -1)y = 0.
\end{gather}
Since $p= \frac{w'}{w} = 0$, we find that $w \in C_K$. Without lost of generality, we can assume that $w=1$. Applying Corollary~\ref{corso3}, we obtain a fundamental matrix of solutions for the rigid solid problem:
\begin{gather*}
\mathbf{Z} = \begin{pmatrix}
y_1^2 - (y_1')^2 & y_1 y_2 - y_1' y_2' & y_2^2 -(y_2')^{2} \\[3pt]
{\rm i} \big(y_1^2 + (y_1')^{2}\big) & {\rm i} \big(y_1 y_2 + y_1' y_2'\big) &
{\rm i} \big(y_2^2 + (y_2')^{2}\big) \\[3pt] -2y_1 y_1' & - y_1 y_2' - y_1' y_2 & -2 y_2 y_2' \end{pmatrix}\!,
\end{gather*}
where $\{ y_1, y_2 \}$ is a basis of solutions of equation~\eqref{eq:RS-eq1}.

Next, we consider the second approach for $\mathfrak{so}(3, C_K)$ systems, given by equation~\eqref{eq:so-system-2}. The identification of matrices $A$ and $\widehat{\Omega}_0$ leads to the degenerate situation: $\omega_1 = - (1/w + wq)$, $\omega_2 = 0$ and $1/w - wq =0$. Thus, $w = -2/ \omega_1$ and $q = \omega_1^2 /4$. This case also produces a coupled situation, this time between $w$ and $q$: $q = 1/w^2$. The fact that $\omega_2 = 0$ means that we are only considering rigid transformations in the line, which is a simplification of the problem. The second-order linear differential equation associated to the system in this case is
\begin{gather}\label{eq:RS-eq2}
 \mathcal{L} y = y'' - \dfrac{\omega_1'}{\omega_1} y' +\dfrac{\omega_1^2}{4}y =0.
\end{gather}
Under this assumptions, a fundamental matrix for this system, given by~\eqref{eq:sol-fund-so3-2}, becomes
\begin{gather*}
\mathbf{Z}_1 = \begin{pmatrix}
 y_1^2 + \dfrac{4 (y_1')^2}{\omega_1^2} & y_1 y_2 + \dfrac{4 y_1' y_2'}{\omega_1^2} & y_2^2 + \dfrac{4 (y_2')^2}{\omega_1^2} \\[10pt]
 - \dfrac{4{\rm i} y_1 y_1'}{\omega_1} & -\dfrac{2{\rm i} (y_1 y_2' + y_1' y_2)}{\omega_1} & - \dfrac{4{\rm i} y_2 y_2'}{\omega_1} \\[10pt]
 {\rm i}\bigg (y_1^2 - \dfrac{4 (y_1')^2}{\omega_1^2} \bigg ) & {\rm i} \bigg (y_1 y_2 - \dfrac{4 y_1' y_2'}{\omega_1^2} \bigg ) & {\rm i}\bigg (y_2^2 - \dfrac{4 (y_2')^2}{\omega_1^2} \bigg )
\end{pmatrix}
\end{gather*}
for $\{ y_1, y_2\}$ a basis of solutions of equation~\eqref{eq:RS-eq2}.

Finally, we consider the following perturbation for the rigid solid system, according to the first approach for $\mathfrak{so}(3, C_K)$ systems (see
	equation~\eqref{eq:system so3}):
\begin{gather*}
\begin{pmatrix} \gamma_1 \\ \gamma_2 \\ \gamma_3 \end{pmatrix}' = - \left (\begin{pmatrix}
0 & 0 & \omega_2 \\ 0 & 0 & - \omega_1 \\ -\omega_2 & \omega_1 & 0 \end{pmatrix} + m \begin{pmatrix} 0 & 0 & -1 \\ 0 & 0 & {\rm i} \\ 1 & -{\rm i} & 0 \end{pmatrix} \right ) \cdot
\begin{pmatrix} \gamma_1 \\ \gamma_2 \\ \gamma_3 \end{pmatrix}\! .
\end{gather*}
Following the philosophy of Darboux in the same vein as we did for the Frenet--Serret system, we can construct an infinite chain of such perturbed Poisson equations for the rigid body problem by applying Proposition~\ref{secthdars031} (again, we restrict ourselves to the case $r=1$). The Darboux transformation is given by
\begin{gather*}
T_{1,m} = \frac{1}{2} \begin{pmatrix}
-{\nu}^{2} + 2 \theta_0 ^{2} -1 & {\rm i} \big({\nu}^{2} -1\big) & 2\theta_0 (1 - \nu) \\
 {\rm i} \big({\nu}^{2} -1\big) &{\nu}^{2}+ 2{\theta_0 }^{2}+1 & 2 {\rm i} \theta_0 (1 +\nu) \\
 2 \theta_0 (\nu -1) &-2 {\rm i} \theta_0 (\nu +1) & 2 \big(\nu + \theta_0 ^2\big)
 \end{pmatrix}
\end{gather*}
for $\theta_0 = \frac{y'}{y}$ and $\nu=m + \theta_0 ^2$, where $y$ is a solution of equation~\eqref{eq:RS-eq1}.
This transformation factors~as
\begin{gather*}
T_{1,m} = \begin{pmatrix}
-m^2 & \theta_0 m & - \theta_0^2 +1 \\ {\rm i}m^2 & -{\rm i} \theta_0 m & {\rm i} + {\rm i} \theta_0^2 \\ 0 & -m & 2\theta_0 \end{pmatrix} \cdot \begin{pmatrix} \frac{1}{2} & - \frac{\rm i}{2} & 0 \\ -\theta_0 & {\rm i} \theta_0 & -1 \\ \frac{1}{2} \big(\theta_0^2 -1\big) & - \frac{\rm i}{2} \big(\theta_0^2 +1\big) & \theta_0
\end{pmatrix}\!.
\end{gather*}

\section{Final remarks}

In this paper, we have elaborated a methodology which shows how, using tensor construction on ${\rm SL}(2, C_K)$, we can explicitly (and algorithmically) obtain Darboux transformations for higher-order linear differential systems such as $\Sym^2({\rm SL}(2, C_K))$-systems or ${\rm SO}(3,C_K)$ systems; the formulas are summarized in the diagrams~\eqref{eq:diag-1} and~\eqref{eq:diag-2}.

Our tool to achieve this is the observation that Darboux transformations can be viewed as gauge transformations and hence may be extended using the tools of Tannakian constructions.\looseness=1

Our approach allows to solve $\Sym^2(\mathrm{SL}(2, C_K))$ systems and produces two Darboux transformations for these kind of systems. In a natural way, we extend this to $\mathfrak{so}(3, C_K)$ systems: we transform an $\mathfrak{so}(3, C_K)$ system into a construction on a second-order linear differential equation by means of the classical isomorphism between the Lie algebras $\mathfrak{so}(3, C_K)$ and $\mathfrak{sl} (2, C_K)$.
This series of transformations give explicit Darboux transformations for $\mathfrak{so}(3, C_K)$ systems, as well as first integrals and simple formulas (already known to experts) for solving them.\looseness=1

These constructions are applied to toy formalisms for supersymmetric quantum mechanics in the non-relativistic case. We have constructed systems-like Schr\"odinger equations following these approaches.
Some well-known $\mathfrak{so}(3, C_K)$ systems such as Frenet--Serret formulas and the rigid solid problem are also included in these constructions.

We notice that the two approaches for $\mathfrak{so}(3, C_K)$ systems developed in this article are not strictly equivalent, since they produce different formulas. This can be seen in Section~\ref{sec:applic-so3}, where two applications to $\mathfrak{so}(3, C_K)$ systems are showed. In the first one, the Frenet--Serret formulas, the first approach leads to a degenerate situation for curves with torsion $\tau= -2i$, whilst the second approach allows us to deal with any curve. However, in the rigid solid application, it is the other way around: the first approach produces a situation that, despite being a coupled case, is more general than the one given by the second approach, which restricts to rigid transformations in the line.

The philosophy developed in this work is based on the Tannakian constructions. Hence, it can straightforwardly allow one to construct Darboux transformations for other differential systems of order higher than two, namely those which can be obtained from a tensor construction on $\mathrm{SL}(2, C_K)$; examples of such systems can be found in~\cite{CoPl94d}.
For any such system, the methodology exposed here allows to construct formulas, solve via solutions of second-order equations and extend Darboux transformations to these families.

\subsection*{Acknowledgements}

The first author thanks the hospitality of XLim and suggestions of J.J.~Morales-Ruiz during the initial stage of this work.
He was supported in the final stage of this paper by the FONDOCYT grants 2022-1D2-90 and 2022-1D2-091 from the Dominican Government (MESCYT).
The third author thanks Autonomous University of Madrid for the financial support for a research stay at XLim, where she started to work in this article. She also thanks the hospitality of XLim and the support of J.J.~Morales-Ruiz to participate in this work.

This work was partially supported by the grant TIN2016-77206-R from the Spanish Government, co-financed by the European Regional Development Fund. The third author received a~postdoctoral grant (PEJD-2018-POST/TIC-9490) from Universidad Nacional de Educaci\'on a~Distancia (UNED), co-financed by the Regional Government of Madrid and the Youth Employment Initiative (YEI) of the European Social Fund.

Authors gratefully acknowledge the referees for
their helpful comments and further references which resulted in an improvement of the preliminary manuscript.

\pdfbookmark[1]{References}{ref}
\LastPageEnding

\end{document}